\documentclass[10pt,a4paper]{amsart}
\usepackage{amssymb,amsmath}
\usepackage[american,french]{babel}
\usepackage{amsopn}
\usepackage{graphicx}
\usepackage{fancyhdr}
\usepackage[T1]{fontenc}
\title{About quadratic differential operators}
\newcommand{\rr}{\mathbb{R}}
\newcommand{\eps}{\varepsilon}
\newcommand{\nn}{\mathbb{N}}
\newcommand{\cc}{\mathbb{C}}

\newcommand{\lde}{L^2(\rr^n)}

\def\init{\setcounter{equa}{0}}
\def\inc{\stepcounter{equa}}
\def\num{\tag{\thesubsection.\theequa}}

\def\wrtext#1{\relax\ifmmode{\leavevmode\hbox{#1}}\else{#1}\fi}
\def\abs#1{\left|#1\right|}
\def\begeq{\begin{equation}}
\def\endeq{\end{equation}}

\def\part#1{\frac{\partial}{\partial #1}}

\def\norm#1{||\,#1\,||}

\newcommand{\real}{\mathbb{R}}
\newcommand{\comp}{\mathbb{C}}

\newcommand{\nat}{\mathbb{N}}

\renewcommand{\Re}{\textrm{Re }}
\renewcommand{\Im}{\textrm{Im }}
\renewcommand{\exp}{\mbox{\rm exp\,}}

\def\Re{{\rm Re\,}}
\def\Im{{\rm Im\,}}

\begin{document}
\newcounter{equa}
\selectlanguage{american}
\begin{center}
{\large\textbf{SPECTRA AND SEMIGROUP SMOOTHING FOR NON-ELLIPTIC QUADRATIC OPERATORS}\\
\bigskip
\medskip
Michael \textsc{Hitrik}, Karel \textsc{Pravda-Starov}\\
\bigskip
University of California, Los Angeles}
\end{center}
\bigskip
\bigskip

\newtheorem{lemma}{Lemma}[subsection]
\newtheorem{definition}{Definition}[subsection]
\newtheorem{proposition}{Proposition}[subsection]
\newtheorem{theorem}{Theorem}[subsection]

\textbf{Abstract.} We study non-elliptic quadratic differential operators. Quadratic differential operators are
non-selfadjoint operators defined in the Weyl quantization by complex-valued quadratic symbols. When the real part of their Weyl symbols is a non-positive quadratic form,
we point out the existence of a particular linear subspace in the phase space intrinsically associated to their
Weyl symbols, called a singular space, such that when the singular space has a symplectic
structure, the associated heat semigroup is smoothing in every direction of its symplectic orthogonal space. When the
Weyl symbol of such an operator
is elliptic on the singular space, this space is always symplectic and we prove that the spectrum of the operator
is discrete and can be described as in the case of
global ellipticity. We also describe the large time behavior of contraction semigroups generated by these operators.

\medskip


\medskip


\section{Introduction}

\subsection{Miscellaneous facts about quadratic differential operators}
\init
Since the classical work by J. Sj\"ostrand~\cite{sjostrand}, the study of spectral properties of quadratic diffe\-rential
operators has played a basic r\^ole in the analysis of partial differential operators with double characteristics.
Roughly speaking, if we have, say, a classical pseudodifferential operator $p(x,\xi)^w$ on $\rr^n$ with the Weyl
symbol $p(x,\xi)=p_m(x,\xi)+p_{m-1}(x,\xi)+\ldots$ of order $m$, and if $X_0=(x_0,\xi_0)\in \rr^{2n}$ is a point
where $p_m(X_0)=dp_m(X_0)=0$ then it is natural to consider the quadratic form $q$ which begins the Taylor expansion of
$p_m$ at $X_0$. The study of a priori estimates for $p(x,\xi)^w$, such as hypoelliptic estimates of the form
$$
\norm{u}_{m-1}\leq C_K \left(\norm{p(x,\xi)^w u}_0+\norm{u}_{m-2}\right),\quad u\in C^{\infty}_0(K), \quad
K\subset \subset \rr^n,
$$
then often depends on the spectral analysis of the quadratic operator $q(x,\xi)^w$.
See also~\cite{hypoelliptic}, as well as Chapter 22 of~\cite{hormander} together with further references given there.
In~\cite{sjostrand}, the spectrum of a general quadratic differential operator has been determined, under the
basic assumption of global ellipticity of the associated quadratic form.

%\medskip
%Recent developments in spectral and pseudospectral theory of non-selfadjoint operators~\cite{DeSjZw},~\cite{ET}
%have brought about a renewed interest in the analysis of quadratic differential operators. In a series of recent
%works~\cite{mz},~\cite{karel}, the second author has analyzed the spectral instability and decay properties of the
%heat semigroup associated to an elliptic quadratic differential operator $q(x,\xi)^w$ such that $\Re q\leq 0$
%does not vanish identically. Here too the precise information concerning the spectrum of $q(x,\xi)^w$ is important
%and has to be supplemented by a study of resolvent estimates.

\medskip
Now there exist many situations where one is naturally led to consider non-self\-ad\-joint quadratic differential
operators whose symbols are not elliptic but rather satisfy certain weaker conditions. An example particularly relevant
to the following discussion is obtained if one considers the Kramers-Fokker-Planck operator with a quadratic
potential~\cite{HN}. The corresponding (complex-valued) symbol is not elliptic, but nevertheless, the operator has
discrete spectrum and the associated heat semigroup is well behaved in the limit of large times --- see~\cite{HeSjSt}.

\medskip
The purpose of the present paper is to provide a proof of a number of fairly general results concerning the
spectral and semigroup properties for the class of quadratic differential operators in the case when the global
ellipticity fails. Specifically, and as alluded to above, we shall consider the class of pseudodifferential operators
defined by the Weyl quantization formula,
\begin{equation}\label{3}\inc
q(x,\xi)^w u(x) =\frac{1}{(2\pi)^n}\int_{\rr^{2n}}{e^{i(x-y).\xi}q\Big(\frac{x+y}{2},\xi\Big)u(y)dyd\xi}, \num
\end{equation}
for some symbols $q(x,\xi)$, where $(x,\xi) \in \rr^{n} \times \rr^n$ and $n \in \nn^*$, which
are complex-valued quadratic forms. Since the symbols are quadratic forms, the corresponding operators in
(\ref{3}) are in fact differential operators. Indeed, the Weyl quantization of the quadratic symbol
$x^{\alpha} \xi^{\beta}$, with $(\alpha,\beta) \in \nn^{2n}$ and $|\alpha+\beta| \leq 2$, is the differential operator
$$
\frac{x^{\alpha}D_x^{\beta}+D_x^{\beta} x^{\alpha}}{2}, \ D_x=i^{-1}\partial_x.
$$
Let us also notice that since the Weyl symbols in (\ref{3}) are complex-valued, the quadratic differential operators
are a priori formally non-selfadjoint.

In this paper, we shall first study the properties of contraction semigroups generated by quadratic
differential operators whose Weyl symbols have a non-positive real part,
\begin{equation}\label{smm1}\inc
\textrm{Re }q \leq 0. \num
\end{equation}
Our first goal is to point out the existence of a linear subspace $S$ in $\rr_x^n \times \rr_{\xi}^n$, which will be
called the singular space and which is defined in terms of the Hamilton map of the Weyl symbol $q$, such that when $S$
has a \textit{symplectic structure}, the associated heat equation
\begin{equation}\label{sm1}\inc
\left\lbrace
\begin{array}{c}
\displaystyle \frac{\partial u}{\partial t}(t,x)-q(x,\xi)^w u(t,x)=0  \\
u(t,\textrm{\textperiodcentered})|_{t=0}=u_0 \in L^2(\rr^n),
\end{array} \right.
 \num
\end{equation}
is smoothing in every direction of the orthogonal complement $S^{\sigma \perp}$ of $S$ with respect to the canonical
symplectic form $\sigma$ on $\rr^{2n}$,
\begin{equation}
\label{11}\inc
\sigma \big{(}(x,\xi),(y,\eta) \big{)}=\xi.y-x.\eta, \ (x,\xi) \in \rr^{2n},  (y,\eta) \in \rr^{2n}. \num
\end{equation}
We shall also describe the large time behavior of contraction semigroups
$$e^{tq(x,\xi)^w}, \ t \geq 0,$$
associated to (\ref{sm1}). When the Weyl symbol $q$ satisfies (\ref{smm1}) and an assumption of
\textit{partial} ellipticity, namely when $q$ is elliptic on the singular space $S$ in the sense that
\begin{equation}\label{sm2}\inc
(x,\xi) \in S, \ q(x,\xi)=0 \Rightarrow (x,\xi)=0, \num
\end{equation}
then $S$ is automatically symplectic, and we prove that the spectrum of the quadratic differential operator
$q(x,\xi)^w$ is only composed of a countable number of eigenvalues of finite multiplicity, with its structure
similar to the one known in the case of global ellipticity~\cite{sjostrand}.

\medskip
\noindent
It seems to us that the singular space $S$ introduced in this paper plays a basic r\^ole
in the understanding of non-elliptic quadratic differential operators. Its study may therefore
be also particularly relevant in the analysis of general pseudodifferential operators with double characteristics,
when the ellipticity of their quadratic approximations fails.

\medskip
Before giving the precise statements of these results, let us begin by recalling some facts and notation about
quadratic differential operators. Let
\begin{eqnarray*}
q : \rr_x^n \times \rr_{\xi}^n &\rightarrow& \cc\\
 (x,\xi) & \mapsto & q(x,\xi),
\end{eqnarray*}
be a complex-valued quadratic form with a non-positive real part,
\begin{equation}\label{inf1}\inc
\textrm{Re }q(x,\xi) \leq 0, \ (x,\xi) \in \rr^{2n}, n \in \nn^*. \num
\end{equation}
We know from \cite{mehler} (p.425) that the
maximal closed realization of  the operator $q(x,\xi)^w$, i.e., the operator on $L^2(\rr^n)$ with the domain
$$\{u \in L^2(\rr^n) :  q(x,\xi)^w u \in L^2(\rr^n)\},$$
coincides with the graph closure of its restriction to $\mathcal{S}(\rr^n)$,
$$q(x,\xi)^w : \mathcal{S}(\rr^n) \rightarrow \mathcal{S}(\rr^n),$$
and that every quadratic differential operator whose Weyl symbol has a non-positive real part,
generates a contraction semigroup. The Mehler formula proved by L. Hörmander in \cite{mehler} gives
an explicit expression for the Weyl symbols of these contraction semigroups.

Associated to the quadratic symbol $q$ is the numerical range $\Sigma(q)$ defined as the closure in the
complex plane of all its values,
\begin{equation}\label{9}\inc
\Sigma(q)=\overline{q(\rr_x^n \times \rr_{\xi}^n)}. \num
\end{equation}
We also recall~\cite{hormander} that the Hamilton map $F \in M_{2n}(\cc)$ associated to the quadratic form $q$
is the map uniquely defined by the identity
\begin{equation}\label{10}\inc
q\big{(}(x,\xi);(y,\eta) \big{)}=\sigma \big{(}(x,\xi),F(y,\eta) \big{)}, \ (x,\xi) \in \rr^{2n},  (y,\eta) \in \rr^{2n}, \num
\end{equation}
where $q\big{(}\textrm{\textperiodcentered};\textrm{\textperiodcentered} \big{)}$ stands for the polarized form
associated to the quadratic form $q$. It follows directly from the definition of the Hamilton map $F$ that
its real part $\textrm{Re } F$ and its imaginary part $\textrm{Im }F$ are the Hamilton maps associated
to the quadratic forms $\textrm{Re } q$ and $\textrm{Im }q$, respectively. Next, (\ref{10}) shows that a
Hamilton map is always skew-symmetric with respect to $\sigma$. This is just a consequence of the
properties of skew-symmetry of the symplectic form and symmetry of the polarized form,
\begin{equation}\label{12}\inc
\forall X,Y \in \rr^{2n}, \ \sigma(X,FY)=q(X;Y)=q(Y;X)=\sigma(Y,FX)=-\sigma(FX,Y).\num
\end{equation}

Let us now consider the elliptic case, i.e., the case of quadratic differential operators whose Weyl symbols are \textit{globally} elliptic
in the sense that
\begin{equation}\label{sm3}\inc
(x,\xi) \in \rr^{2n}, \ q(x,\xi)=0 \Rightarrow (x,\xi)=0. \num
\end{equation}
In this case, the numerical range of a quadratic form can only take very particular shapes. J. Sjöstrand proved in \cite{sjostrand} (Lemma 3.1) that
if $q$ is a complex-valued elliptic quadratic form on $\rr^{2n}$, with $n \geq 2$, then there exists $z \in \cc^*$
such that $\textrm{Re}(z q)$ is a positive definite quadratic form. If $n=1$, the same result is fulfilled if we
assume besides that $\Sigma(q) \neq \cc$. This shows that the numerical range of an elliptic quadratic form can only
take two shapes. The first possible shape is when $\Sigma(q)$ is equal to the whole complex plane. This case can
only occur in dimension $n=1$. The second possible shape is when $\Sigma(q)$ is equal to a closed angular sector
with a vertex in $0$ and an aperture strictly less than $\pi$ (see \cite{mz} for more details).

We also know that elliptic quadratic differential operators define Fredholm operators
(see Lemma~3.1 in \cite{hypoelliptic} or Theorem 3.5 in \cite{sjostrand}),
\begin{equation}\label{14}\inc
q(x,\xi)^w +z : B \rightarrow \lde, \num
\end{equation}
where $B$ is the Hilbert space
\inc
\begin{align*}\label{14.1}
B= & \ \big\{ u \in \lde : q(x,\xi)^wu \in \lde \big\} \num \\
= & \ \big\{ u \in \lde : x^{\alpha} D_x^{\beta} u \in \lde \ \textrm{if} \ |\alpha+\beta| \leq 2\big\},
\end{align*}
with the norm
$$\|u\|_B^2=\sum_{|\alpha+\beta| \leq 2}{\|x^{\alpha} D_x^{\beta} u\|_{\lde}^2}.$$
Moreover, the index of the operator (\ref{14}) is independent of $z$ and is equal to~$0$ when $n \geq 2$.
In the case $n=1$, the index can take the values $-2$, $0$ or $2$. It vanishes as soon as $\Sigma(q) \neq \cc$.

When $\Sigma(q) \neq \cc$, J.~Sjöstrand has proved in Theorem 3.5 of \cite{sjostrand} (see also Lemma~$3.2$ and Theorem~3.3
in~\cite{hypoelliptic}) that the spectrum of an elliptic quadratic differential operator
$$q(x,\xi)^w : B \rightarrow \lde,$$
is only composed of eigenvalues with finite multiplicity,
\begin{equation}\label{15}\inc
\sigma\big{(}q(x,\xi)^w\big{)}=\Big\{ \sum_{\substack{\lambda \in \sigma(F), \\  -i \lambda \in \Sigma(q)\setminus \{0\}} }
{\big{(}r_{\lambda}+2 k_{\lambda}
\big{)}(-i\lambda) : k_{\lambda} \in \nn}
\Big\}, \num
\end{equation}
where $F$ is the Hamilton map associated to the quadratic form $q$ and $r_{\lambda}$ is the dimension of the space of generalized eigenvectors of $F$ in $\cc^{2n}$
belonging to the eigenvalue $\lambda \in \cc$.

Let us also recall the result proved in \cite{mz} about contraction semigroups generated by elliptic quadratic differential operators whose Weyl symbols have a non-positive real part.
This result shows that, as soon as the real part of their Weyl symbols is a non-zero quadratic form, the norm of contraction semigroups generated by these operators decays
exponentially in time.

In this paper, we study the case when the ellipticity fails. Our second result (Theorem~\ref{theorem2}) extends
the description of the spectra (\ref{15}) to the case of quadratic differential operators whose Weyl symbols
are partially elliptic, but not necessarily globally so. To get this result, we only require that these
symbols have a non-positive real part and are elliptic on their associated singular spaces. We also prove a
result on the exponential decay in time for the norm of contraction semigroups generated by non-elliptic quadratic
differential operators.

Let us now define this singular space. The \textit{singular} space $S$ associated to the symbol $q$ is defined as
the following intersection of the kernels,
\begin{equation}\label{h1}\inc
S=\Big(\bigcap_{j=0}^{+\infty}\textrm{Ker}\big[\textrm{Re }F(\textrm{Im }F)^j \big]\Big) \cap \rr^{2n}, \num
\end{equation}
where the notation $\textrm{Re } F$ and $\textrm{Im }F$ stands respectively for the real part and the imaginary part of the Hamilton map associated to $q$.
Notice that the Cayley-Hamilton theorem applied to $\textrm{Im }F$ shows that
$$(\textrm{Im }F)^k X \in \textrm{Vect}\big(X,...,(\textrm{Im }F)^{2n-1}X\big), \ X \in \rr^{2n}, \ k \in \nn,$$
where $\textrm{Vect}\big(X,...,(\textrm{Im }F)^{2n-1}X\big)$ is the vector space spanned by the vectors $X$, ...,
$(\textrm{Im }F)^{2n-1}X$, and therefore
the singular space is actually equal to the following finite intersection of the kernels,
\begin{equation}\label{h1bis}\inc
S=\Big(\bigcap_{j=0}^{2n-1}\textrm{Ker}\big[\textrm{Re }F(\textrm{Im }F)^j \big]\Big) \cap \rr^{2n}. \num
\end{equation}
The subspace $S$ obviously satisfies the two following properties,
\begin{equation}\label{h3}\inc
(\textrm{Re }F) S=\{0\} \textrm{ and } (\textrm{Im }F)S \subset S. \num
\end{equation}
We can now give the statements of the main results contained in this paper.

\subsection{Statement of the main results}
\init

In the following statements, we consider a complex-valued quadratic form
$$q : \rr_x^n \times \rr_{\xi}^n \rightarrow \cc,$$
with a non-positive real part,
\begin{equation}\label{sm4}\inc
\textrm{Re }q(x,\xi) \leq 0, \ (x,\xi) \in \rr^{2n}, \ n \in \nn^*, \num
\end{equation}
and we denote by $S$ the singular space defined in (\ref{h1}) or (\ref{h1bis}).

Our first result states that when the singular space $S$ has a symplectic structure, in the sense that the restriction of
$\sigma$ to $S$ is nondegenerate, the heat equation (\ref{sm1}) associated to the operator $q(x,\xi)^w$
is smoothing in every direction of its orthogonal complement $S^{\sigma \perp}$
with respect to the canonical symplectic form in $\rr^{2n}$.

\bigskip

\begin{theorem}\label{theorem1}
Let us assume that the singular space $S$ has a symplectic structure.
If $(x',\xi')$ are some linear symplectic coordinates on the symplectic space $S^{\sigma \perp}$,
then for all $t>0$, $N \in \nn$ and $u \in L^2(\rr^n)$,
\begin{equation}\label{sm4bis}\inc
\big((1+|x'|^2+|\xi'|^2)^N\big)^we^{t q(x,\xi)^w}u \in L^2(\rr^n). \num
\end{equation}
\end{theorem}

\bigskip

Let us mention that the assumption about the symplectic structure of $S$ is always fulfilled by any
quadratic symbol $q$ elliptic on $S$, i.e.,
\begin{equation}\label{sm5}\inc
(x,\xi) \in S, \ q(x,\xi)=0 \Rightarrow (x,\xi)=0. \num
\end{equation}
This assumption is therefore always fulfilled for elliptic quadratic differential operators. We will
see that it is also the case for instance for the Kramers-Fokker-Planck operator with a quadratic potential,
which is a non-elliptic operator.

When $q$ is a complex-valued quadratic form with a non-positive real part verifying (\ref{sm5}), we can
give another description of the singular space in terms of the eigenspaces of $F$
associated to its real eigenvalues. Under these assumptions, the set of real eigenvalues of the Hamilton map $F$ can be written as
$$\sigma(F) \cap \rr =\{\lambda_1,...,\lambda_r,-\lambda_1,...,-\lambda_r\},$$
with $\lambda_j \neq 0$ and $\lambda_j \neq \pm \lambda_k$  if $j \neq k$.
The singular space is then the direct sum of the symplectically orthogonal spaces
\begin{equation}\label{sm5bis}\inc
S=S_{\lambda_1} \oplus^{\sigma \perp} S_{\lambda_2} \oplus^{\sigma \perp}... \oplus^{\sigma \perp} S_{\lambda_r}, \num
\end{equation}
where $S_{\lambda_{j}}$ is the symplectic space
\begin{equation}\label{sm5bis1}\inc
S_{\lambda_j}=\big(\textrm{Ker}(F -\lambda_j) \oplus \textrm{Ker}(F+\lambda_j) \big) \cap \rr^{2n}. \num
\end{equation}
These facts will be proved in section \ref{examples}.

Our second result deals with the structure of the spectra for non-elliptic quadratic differential operators. This
result extends the description of the spectra (\ref{15}) proved by J.~Sjöstrand in \cite{sjostrand} (Theorem 3.5)
for elliptic quadratic differential operators to the case of quadratic differential operators which are only
partially elliptic. To get this description, we only require in addition to the assumption (\ref{sm4})
the property of partial ellipticity (\ref{sm5}) for their Weyl symbols.

\bigskip

\begin{theorem}\label{theorem2}
If $q$ is a complex-valued quadratic form with a non-positive real part and if $q$ is elliptic on $S$,
$$(x,\xi) \in S, \ q(x,\xi)=0 \Rightarrow (x,\xi)=0,$$
then the spectrum of the quadratic differential operator  $q(x,\xi)^w$
is only composed of eigenvalues of finite multiplicity,
\begin{equation}\label{sm6}\inc
\sigma\big{(}q(x,\xi)^w\big{)}=\Big\{ \sum_{\substack{\lambda \in \sigma(F), \\  -i \lambda \in \cc_-
\cup (\Sigma(q|_S) \setminus \{0\})
} }
{\big{(}r_{\lambda}+2 k_{\lambda}
\big{)}(-i\lambda) : k_{\lambda} \in \nn}
\Big\}, \num
\end{equation}
where $F$ is the Hamilton map associated to the quadratic form $q$, $r_{\lambda}$ is the dimension of the space of generalized eigenvectors of $F$ in $\cc^{2n}$
belonging to the eigenvalue $\lambda \in \cc$,
$$\Sigma(q|_S)=\overline{q(S)} \textrm{ and } \cc_-=\{z \in \cc : \emph{\textrm{Re }}z<0\}.$$
\end{theorem}

\bigskip

Since the singular space $S$ is distinct from the whole phase space as soon as the real part of $q$ is not
identically equal to zero, Theorem \ref{theorem2} is a generalization of the result
proved by J. Sjöstrand for elliptic quadratic differential operators.

Finally, we give a result concerning the large time behavior of contraction semigroups generated by
non-elliptic quadratic differential operators, which extends the result obtained by the second author in \cite{mz}.

\bigskip

\begin{theorem}\label{theorem3}
Let us consider a complex-valued quadratic form
$$q : \rr_x^n \times \rr_{\xi}^n \rightarrow \cc, \ n \in \nn^*,$$
with a non-positive real part, such that its singular space $S$ has a symplectic structure. Then, the following assertions
are equivalent:
\begin{itemize}
\item[$(i)$] The norm of the contraction semigroup generated by the operator $q(x,\xi)^w$ decays exponentially in time,
$$\exists M>0, \exists a >0, \forall t \geq 0, \ \|e^{tq(x,\xi)^w}\|_{\mathcal{L}(L^2)} \leq M e^{-a t}.$$
\item[$(ii)$] The real part of the symbol $q$ is a non-zero quadratic form
$$\exists (x_0,\xi_0) \in \rr^{2n}, \ \emph{\textrm{Re }}q(x_0,\xi_0) \neq 0.$$
\item[$(iii)$] The singular space is distinct from the whole phase space $S \neq \rr^{2n}$.
\end{itemize}
\end{theorem}

\bigskip

Since the assumption about the symplectic structure of the singular space $S$ is always fulfilled when the
symbol $q$ verifies (\ref{sm5}), Theorem \ref{theorem3} is
a generalization of the result proved in \cite{mz} for elliptic quadratic differential operators.

Let us also notice that we cannot drop completely the assumption about the symplectic structure of the singular space.  Indeed, let us consider
the quadratic differential operator defined in the Weyl quantization by the symbol
$$q(x,\xi)=-x^2.$$
This operator is just the operator of multiplication by $-x^2$, which generates the contraction semigroup
$$e^{t q(x,\xi)^w}u=e^{-t x^2}u, \ t \geq 0, \ u \in L^2(\rr^n),$$
whose norm is identically equal to 1,
$$\|e^{t q(x,\xi)^w}\|_{\mathcal{L}(L^2)}=1, \ t \geq 0.$$

\bigskip

\noindent
\textit{Remark.} Let us mention that our proof will show in particular that when $q$ is a complex-valued quadratic form
on $\rr^{2n}$, $n \geq 1$, with a non-positive real part and a zero
singular space $S=\{0\}$, then
$$e^{t q(x,\xi)^w}=e^{t q(x,\xi)^w} \Pi_a+\mathcal{O}_a(e^{-at}), \ t \geq 0,$$
in the space $\mathcal{L}(L^2)$ of bounded operators on $L^2(\rr^n)$, for any $a>0$ such that
\begin{align*}
& \ \sigma\big(q(x,\xi)^w\big) \cap \{ z \in \cc :  \Re z =-a\} \\
= & \ \Big\{ \sum_{\substack{\lambda \in \sigma(F), \\  \textrm{Re}(-i \lambda) <0
} }
{\big{(}r_{\lambda}+2 k_{\lambda}
\big{)}(-i\lambda) : k_{\lambda} \in \nn}
\Big\}\cap \{ z \in \cc :  \Re z=-a \}=\emptyset,
\end{align*}
where $\Pi_a$ stands for the finite rank spectral projection associated to the following eigenvalues of the operator $q(x,\xi)^w$,
\begin{align*}
& \ \sigma\big(q(x,\xi)^w\big) \cap \{ z \in \cc :  -a \leq \Re z \} \\
= & \ \Big\{ \sum_{\substack{\lambda \in \sigma(F), \\  \textrm{Re}(-i \lambda) <0
} }
{\big{(}r_{\lambda}+2 k_{\lambda}
\big{)}(-i\lambda) : k_{\lambda} \in \nn}
\Big\}\cap \{ z \in \cc :  -a \leq \Re z \},
\end{align*}
where $F$ is the Hamilton map associated to the quadratic form $q$ and $r_{\lambda}$ is the dimension of the space of generalized eigenvectors of $F$ in $\cc^{2n}$
belonging to the eigenvalue $\lambda \in \cc$.

\bigskip

Let us now explain the key arguments in our proofs of these theorems.

\subsection{Structure of the proof}
\init

Our main assumption about the symplectic structure of the singular space $S$ fulfilled in the assumptions of all the three theorems allows us to
find some symplectic coordinates $(x',\xi')$ in $S^{\sigma \perp}$ and $(x'',\xi'')$ in
$S$ such that the complex-valued quadratic form $q$ verifying (\ref{sm4}) can be written as the sum of two quadratic forms with a tensorization of the variables $(x',\xi')$ and $(x'',\xi'')$,
$$q=q|_S+q|_{S^{\sigma \perp}}, \ (x,\xi)=(x',x'';\xi',\xi'') \in \rr^{2n},$$
where the first quadratic form $q|_S$ is equal to
$$q|_S=i \tilde{q}|_S,$$
with $\tilde{q}|_S$ a real-valued quadratic form; and where the second quadratic form $q|_{S^{\sigma \perp}}$ is a complex-valued quadratic form with a non-positive real part.
This real part is not in general negative definite (unless the real part of $q$ is). However, it follows from
the definition of the singular space $S$ that the average of the real part of the quadratic form
$q|_{S^{\sigma \perp}}$ by the flow generated by the Hamilton vector field of its imaginary part,
$H_{\textrm{Im}q|_{S^{\sigma \perp}}}$,
$$
\langle \textrm{Re }q|_{S^{\sigma \perp}}\rangle_T(X')=\frac{1}{2T}\int_{-T}^T{\textrm{Re }q|_{S^{\sigma \perp}}(e^{tH_{\textrm{Im}q|_{S^{\sigma \perp}}}}X')dt}, \ T>0, \ X'=(x',\xi'),$$
is negative definite. Studying the contraction semigroup
\begin{equation}
\label{tom1}\inc
e^{tq|_{S^{\sigma \perp}}^w}, \ t \geq 0,
\num
\end{equation}
generated by the operator $q|_{S^{\sigma \perp}}^w$, on the FBI-Bargmann transform side, we prove, using
the averaging property just mentioned, that (\ref{tom1}) is compact
and strongly regularizing for every $t >0$. This compactness result is really the key point in our proofs of
the three theorems, and their complete statements then follow from a small additional amount of work.

\subsection{Some examples}\label{examples}
\init

In this section, we prove that if a quadratic symbol $q$ is elliptic on its singular space then the
singular space always has a symplectic structure. We also check that this
property is fulfilled for the Kramers-Fokker-Plank operator with a quadratic potential.

\subsubsection{Partially elliptic quadratic differential operators}

Let us consider the case of quadratic differential operators whose Weyl symbols are elliptic on their singular spaces. Let
$$q : \rr_x^n \times \rr_{\xi}^n \rightarrow \cc, \ n \in \nn^*,$$
be a complex-valued quadratic form, which is elliptic on its singular space $S$,
\begin{equation}\label{sm12bis}\inc
(x,\xi) \in S, \ q(x,\xi)=0 \Rightarrow (x,\xi)=0. \num
\end{equation}
We want to prove that
$$S=\Big(\bigcap_{j=0}^{2n-1}\textrm{Ker}\big[\textrm{Re }F(\textrm{Im }F)^j \big]\Big) \cap \rr^{2n},$$
has a symplectic structure. This fact follows from some arguments similar to those used in \cite{mz} (Lemma 3).

We can assume that $S \neq \{0\}$ since the space $\{0\}$ is obviously symplectic. Let us therefore consider $X_0 \in S \setminus \{0\}$.
We define
\begin{equation}\label{2.3.101}\inc
\left\lbrace
\begin{array}{c}
e_1=X_0\\
\displaystyle \eps_1=-\frac{1}{\textrm{Im } q(X_0)}\textrm{Im } F X_0.
 \end{array} \right. \num
\end{equation}
This is possible, since from (\ref{10}) and (\ref{h3}), we have
$$\textrm{Re } q(X_0)=\sigma(X_0,\textrm{Re }F X_0)=0,$$
and the ellipticity of $q$ on $S$ implies that
$$\textrm{Im } q(X_0) \neq 0,$$
as $X_0 \in S \setminus \{0\}$. By using the skew-symmetry of the Hamilton map $\textrm{Im }F$ (see (\ref{12})), it follows that
$$\sigma(\eps_1,e_1)=\sigma\big{(}- (\textrm{Im } q(X_0))^{-1} \textrm{Im } F X_0,X_0\big{)}=(\textrm{Im } q(X_0))^{-1}\sigma(X_0,\textrm{Im } F X_0)=1,$$
which shows that the system $(e_1,\eps_1)$ is symplectic. We get from (\ref{h3}) and (\ref{2.3.101}) that
$$\textrm{Vect}(e_1,\eps_1) \subset S.$$
If $S=\textrm{Vect}(e_1,\eps_1)$, the singular space $S$ is symplectic. If it is not the case, so that,
$$S \neq \textrm{Vect}(e_1,\eps_1),$$
we can continue our construction of a symplectic basis for $S$ by considering
$$X_1 \in S \setminus \textrm{Vect}(e_1,\eps_1)$$
and
\begin{equation}\label{hh3}\inc
\tilde{X}_1=X_1+\sigma(X_1,\eps_1)e_1-\sigma(X_1,e_1)\eps_1 \in S \setminus \textrm{Vect}(e_1,\eps_1). \num
\end{equation}
Let us set
\begin{equation}\label{2.3.107}\inc
\left\lbrace
\begin{array}{c}
e_2=\tilde{X}_1\\
\displaystyle \eps_2=-\frac{1}{\textrm{Im } q(\tilde{X}_1)}\big( \textrm{Im } F \tilde{X}_1+\sigma(\textrm{Im }F \tilde{X}_1,\eps_1)e_1-\sigma(\textrm{Im }F \tilde{X}_1,e_1)\eps_1\big),
\end{array} \right. \num
\end{equation}
which is again possible according to (\ref{h3}) and the assumption of ellipticity on $S$ since
$$\textrm{Re } q(\tilde{X}_1)=\sigma(\tilde{X}_1,\textrm{Re }F \tilde{X}_1)=0,$$
because $\tilde{X}_1 \in S \setminus \{0\}$. Then, we can directly verify by using (\ref{hh3}) and (\ref{2.3.107})
that $(e_1,e_2,\eps_1,\eps_2)$ is a symplectic system. By using (\ref{h3}) again, we get
$$\textrm{Vect}(e_1,e_2,\eps_1,\eps_2) \subset S.$$
If $S=\textrm{Vect}(e_1,e_2,\eps_1,\eps_2)$, then $S$ is symplectic. If it is not the case, then
$$S \neq \textrm{Vect}(e_1,e_2,\eps_1,\eps_2),$$
we can again iterate the preceding construction. After a finite number of such iterations,
we obtain with this process a symplectic basis of $S$, proving its symplectic structure.

Let us now consider a complex-valued quadratic form
$$q : \rr_x^n \times \rr_{\xi}^n \rightarrow \cc, \ n \in \nn^*,$$
with a non-positive real part
$$\textrm{Re }q \leq 0,$$
such that (\ref{sm12bis}) is fulfilled and denote by $F$ its Hamilton map. We know from Proposition 4.4 in
\cite{mehler} that the kernel
$\textrm{Ker}(F+\lambda)$ is the complex conjugate of the kernel $\textrm{Ker}(F-\lambda)$ for every $\lambda \in \rr$,
and that the spaces
$$\textrm{Ker}(F-\lambda) \oplus  \textrm{Ker}(F+\lambda),$$
where $\lambda \in \rr^*$, and $\textrm{Ker }F$, are the complexifications of their intersections with $\rr^{2n}$.

Let us set
\begin{equation}\label{venice1}\inc
S_0=(\textrm{Ker }F) \cap \rr^{2n} \num
\end{equation}
and
\begin{equation}\label{venice2}\inc
S_{\lambda}=\big(\textrm{Ker}(F-\lambda) \oplus  \textrm{Ker}(F+\lambda)\big) \cap \rr^{2n}, \num
\end{equation}
for $\lambda \in \rr^*$. Proposition 4.4 in \cite{mehler} also shows that
$$\textrm{Re }  F \ \textrm{Ker}(F \pm \lambda)=\{0\},$$
for all $\lambda \in \rr$. This implies that
\begin{equation}\label{venice3}\inc
(\textrm{Re } F)S_{\lambda}=\{0\} \textrm{ and } (\textrm{Im } F)S_{\lambda} \subset S_{\lambda}, \num
\end{equation}
and proves in view of (\ref{h1bis}) that for all $\lambda \in \rr$,
\begin{equation}\label{venice4}\inc
S_{\lambda} \subset S. \num
\end{equation}
If $0 \in \sigma(F) \cap \rr$, this would imply that $S_0 \neq 0$. Since from (\ref{venice1}),
$$q(X)=\sigma(X,FX)=0,$$
for all $X \in S_0$, the inclusion (\ref{venice4}) would then contradict our assumption of ellipticity on the singular space (\ref{sm12bis}). This proves that the set of real eigenvalues of the
Hamilton map $F$ can be written as
\begin{equation}\label{venice5}\inc
\sigma(F) \cap \rr =\{\lambda_1,...,\lambda_r,-\lambda_1,...,-\lambda_r\}, \num
\end{equation}
with $\lambda_j \neq 0$ and $\lambda_j \neq \pm \lambda_k$ if $j \neq k$.

Let us now check that the spaces $S_{\lambda_j}$, $j=1,...,r$, are symplectic. Let $X_0$ be in
$S_{\lambda_j}$ such that for all $Y \in S_{\lambda_j}$,
$$\sigma(X_0,Y)=0.$$
It follows that for all $Y$ and  $Z$ in $S_{\lambda_j}$,
$$\sigma(X_0,Y+iZ)=0,$$
which induces that
$$\forall X \in \textrm{Ker}(F-\lambda_j) \oplus  \textrm{Ker}(F+\lambda_j), \ \sigma(X_0,X)=0,$$
because $ \textrm{Ker}(F-\lambda_j) \oplus  \textrm{Ker}(F+\lambda_j)$ is a complexification of $S_{\lambda_j}$.
On the other hand, since $X_0 \in S_{\lambda_j}$, we have $F X_0 \in \textrm{Ker}(F-\lambda_j) \oplus  \textrm{Ker}(F+\lambda_j)$, which implies that
$$q(X_0)=\sigma(X_0,F X_0)=0.$$
We then deduce from the ellipticity of $q$ on the singular space (\ref{sm12bis}) and (\ref{venice4}) that $X_0=0$, which proves the symplectic structure of the space $S_{\lambda_j}$.

Let us now assume that there exists another real eigenvalue $\lambda_k$ of $F$ distinct from $\lambda_j$ and $-\lambda_j$. We already know that this eigenvalue $\lambda_k$
is necessarily non-zero. Let
$$X \in \textrm{Ker}(F-\eps_1 \lambda_j) \textrm{ and } Y \in \textrm{Ker}(F-\eps_2 \lambda_k),$$
with $\eps_1, \eps_2 \in \{\pm 1\}$, we obtain from the skew-symmetry property of the Hamilton map $F$ with respect to $\sigma$ that
$$\sigma(X,Y)=\sigma(X,\eps_2^{-1} \lambda_k^{-1} FY)=-\frac{1}{\eps_2 \lambda_k}\sigma(FX,Y)=-\frac{\eps_1 \lambda_j}{\eps_2 \lambda_k} \sigma(X,Y).$$
Since
$$\left|\frac{\eps_1 \lambda_j}{\eps_2 \lambda_k} \right| \neq 1,$$
because $\lambda_j$ and $\lambda_k$ are real numbers such that $\lambda_k \not\in \{\lambda_j,-\lambda_j\}$, we finally deduce that
$$\sigma(X,Y) =0,$$
which proves that the spaces $S_{\lambda_j}$ and $S_{\lambda_k}$ are symplectically orthogonal, and we get from (\ref{venice4}) and (\ref{venice5}) that
\begin{equation}\label{venice6}\inc
S_{\lambda_1} \oplus^{\sigma \perp} S_{\lambda_2} \oplus^{\sigma \perp} ... \oplus^{\sigma \perp} S_{\lambda_r}  \subset S. \num
\end{equation}

Let us prove that the singular space is actually exactly equal to this direct sum of symplectic spaces. We recall that from (\ref{h3}),
$$(\textrm{Re }F) S=\{0\} \textrm{ and } (\textrm{Im }F) S \subset S.$$
Since
$$q(X)=\sigma(X,FX)=i \sigma(X,\textrm{Im }FX), \ X \in S,$$
we deduce from (\ref{sm12bis}) and the lemma 18.6.4 in \cite{hormander} that we can find new symplectic
basis $(\tilde{e}_1,...,\tilde{e}_m,\tilde{\eps}_1,...,\tilde{\eps}_m)$
in the symplectic space $S$ such that
\begin{equation}\label{venice6bis}\inc
q(X)=i \eps \sum_{j=1}^m{\mu_j(\tilde{\xi}_j^2+\tilde{x}_j^2)}, \ X=\tilde{x}_1 \tilde{e}_1+...+ \tilde{x}_m \tilde{e}_m+\tilde{\xi}_1 \tilde{\eps}_1+...+\tilde{\xi}_m \tilde{\eps}_m, \num
\end{equation}
where $\eps \in \{\pm 1\}$ and $\mu_j>0$ for all $j=1,...,m$. Indeed, this is linked to the fact that a real-valued elliptic quadratic form must be positive definite or negative definite.
By computing $F$ from (\ref{venice6bis}), we get that
\begin{equation}\label{venice7}\inc
FX_j=-\eps \mu_j X_j \textrm{ and } F \tilde{X}_j= \eps \mu_j \tilde{X}_j, \num
\end{equation}
if $X_j=\tilde{e}_j+i \tilde{\eps}_j$ and $\tilde{X}_j=\tilde{e}_j-i\tilde{\eps}_j$ for all $j=1,...,m$. The
identities (\ref{venice7}) prove that the singular space is actually equal to
the direct sum of the symplectic spaces $S_{\lambda_j}$ defined in (\ref{venice2}),
$$S=S_{\lambda_1} \oplus^{\sigma \perp} S_{\lambda_2} \oplus^{\sigma \perp} ... \oplus^{\sigma \perp} S_{\lambda_r}.$$

\subsubsection{Kramers-Fokker-Planck operator with a quadratic potential}\label{planck}
Let us consider the Kramers-Fokker-Planck operator~\cite{HN},
$$K=-\Delta_v+\frac{v^2}{4}-\frac{1}{2}+v.\partial_x-\big(\partial_xV(x)\big).\partial_v, \ (x,v) \in \rr^{2},$$
with a quadratic potential
$$V(x)=\frac{1}{2}ax^2, \ a \in \rr^*.$$
We can write
$$K=-q(x,v,\xi,\eta)^w-\frac{1}{2},$$
with
$$q(x,v,\xi,\eta)=-\eta^2-\frac{1}{4}v^2-i(v \xi-a x \eta).$$
This symbol $q$ is a non-elliptic complex-valued quadratic form with a non-positive real part and a numerical range equal to the half-plane
$$\Sigma(q)=\{z \in \cc : \textrm{Re }z \leq 0\}.$$
We can directly check that its Hamilton map $F$ for which
$$q(x,v,\xi,\eta)=\sigma\big((x,v,\xi,\eta),F(x,v,\xi,\eta) \big),$$
is given by
$$F= \left( \begin{array}{cccc}
  0 & -\frac{1}{2}i & 0 & 0 \\
  \frac{1}{2}ai& 0 & 0 & -1 \\
0 & 0 & 0 & - \frac{1}{2}ai \\
0 & \frac{1}{4} & \frac{1}{2}i & 0
  \end{array}
\right),$$
and that the singular space
$$S= \Big(\bigcap_{j=0}^{3}\textrm{Ker}\big[\textrm{Re }F(\textrm{Im }F)^j \big]\Big) \cap \rr^{4},$$
is reduced to the trivial symplectic space $\{0\}$.

\medskip
\noindent
{\bf Acknowledgment}. The research of the first author is supported in part by
the National Science Foundation under grant DMS--0653275 and the Alfred P. Sloan Research Fellowship. He would also like
to thank Joe Viola for a stimulating discussion.

\section{Symplectic decomposition of the symbol}\label{sympl}
\init

In this section, we explain how the main assumption about the symplectic structure of the singular space $S$
fulfilled in the hypotheses of all our three theorems
allows us to tensor the variables in the
symbol $q$ by writing it as a sum of two quadratic forms where the first one is purely imaginary-valued and where
the second one verifies the averaging property of its real part by the flow defined
by the Hamilton vector field of its imaginary part.

Let us consider a complex-valued quadratic form $q$ verifying (\ref{sm4}) and let us assume that the singular space $S$ defined in (\ref{h1bis}) is symplectic. Let us recall
that it is well the case when $q$ verifies (\ref{sm5}). Then, we can find
$\chi$, a real linear symplectic transformation of $\rr^{2n}$, such that
\begin{equation}\label{sm12}\inc
(q \circ \chi)(x,\xi)= q_1(x',\xi')+ iq_2(x'',\xi''), \ (x,\xi)=(x',x'';\xi',\xi'') \in \rr^{2n}, \num
\end{equation}
where $q_1$ is a complex-valued quadratic form on $\rr^{2n'}$ with a non-positive real part
\begin{equation}\label{sm13}\inc
\textrm{Re }q_1 \leq 0, \num
\end{equation}
and $q_2$ is a real-valued quadratic form verifying the following properties:

\medskip

\begin{proposition}\label{proposition1}
The two quadratic forms $q_1$ and $q_2$ satisfy the following properties:
\begin{itemize}
\item[$(i)$] For all $T>0$, the average of the real part of the quadratic form $q_1$ by the flow defined by the Hamilton vector field of $\emph{\textrm{Im }}q_1$,
$$\langle\emph{\textrm{Re }}q_1\rangle_{T}(X')=\frac{1}{2T}\int_{-T}^{T}{\emph{\textrm{Re }}q_1(e^{tH_{\emph{\textrm{Im}}q_1}}X')dt}, \ X'=(x',\xi') \in \rr^{2n'},$$
is negative definite.
\item[$(ii)$] The quadratic form
$$\sum_{j=0}^{2n-1}{\emph{\textrm{Re }} q_1\big((\emph{\textrm{Im }} F_1)^j X'\big)}, \ X'=(x',\xi') \in \rr^{2n'}, $$
where $F_1$ stands for the Hamilton map of $q_1$, is negative definite.
\item[$(iii)$] If the symbol $q$ fulfills an additional assumption of ellipticity on $S$,
\begin{equation}\label{la}\inc
(x,\xi) \in S, \ q(x,\xi)=0 \Rightarrow (x,\xi)=0,\num
\end{equation}
then we can assume that
$$q_2(x'',\xi'')=\eps \sum_{j=1}^{n''}{\lambda_j(\xi_j''^2+x_j''^2)},$$
where $\eps \in \{\pm 1\}$ and $\lambda_j>0$ for all $j=1,...,n''$.
\end{itemize}

\end{proposition}

\medskip

To prove these results, we begin by considering $S^{\sigma \perp}$, the orthogonal complement of
$S$ in $\rr^{2n}$ with respect to the symplectic form and $F$ the Hamilton map of $q$.
The space $S^{\sigma \perp}$ is symplectic because it is the case for $S$.
Moreover, since according to (\ref{h3}), $S$ is stable by the maps $\textrm{Re } F$ and $\textrm{Im } F$,
its orthogonal complement also fulfills these properties.
Indeed, let $X$ be in $S^{\sigma \perp}$. By using (\ref{h3}) and the skew-symmetry of any Hamilton
map with respect to $\sigma$, we get for all $Y \in S$,
$$\sigma(Y,\textrm{Re } F X)+i \sigma(Y,\textrm{Im } F X)=-\sigma(\textrm{Re } F Y,X)-i \sigma(\textrm{Im } F Y,X)=0,$$
because $(\textrm{Re } F) Y \in S$ and $(\textrm{Im } F) Y \in S$.
This induces that for all $Y \in S$,
$$\sigma(Y,\textrm{Re } F X)=\sigma(Y,\textrm{Im } F X)=0,$$
and proves that $(\textrm{Re } F) X \in S^{\sigma \perp  }$ and $(\textrm{Im } F) X \in S^{\sigma \perp  }$.

We can then write the phase space $\rr^{2n}$ as a direct sum of two symplectically orthogonal real symplectic spaces stable by the maps
$\textrm{Re } F$ and $\textrm{Im } F$,
\begin{equation}\label{2.3.89}\inc
\rr^{2n}= S_1 \oplus^{\sigma \perp} S_2, \  (\textrm{Re } F)S_j \subset S_j, \ (\textrm{Im } F)S_j \subset S_j, \num
\end{equation}
for $j \in \{1,2\}$ with
\begin{equation}\label{h5}\inc
S_1=S^{\sigma \perp} \textrm{ and } S_2= S. \num
\end{equation}
Let us now consider a symplectic basis $(e_{1,j}, ...,e_{N_j,j},\eps_{1,j},...,\eps_{N_j,j})$ of $S_j$. By collecting
these two bases, we get a
symplectic basis of $\rr^{2n}$, which allows by using the stability and the orthogonality properties of the spaces $S_j$ to obtain the following decomposition of $q$,
\begin{align*}
 q(x,\xi) = & \ \sigma\Big{(} \sum_{\substack{1 \leq j \leq 2, \\ 1 \leq k \leq N_j}}{(x_{k,j} e_{k,j}+\xi_{k,j} \eps_{k,j})}, F\big{(}
\sum_{\substack{1 \leq j \leq 2,\\  1 \leq k \leq N_j}}{(x_{k,j} e_{k,j}+\xi_{k,j} \eps_{k,j})}\big)\Big)\\
= & \
\sum_{1 \leq j \leq 2}\sigma \Big{(}\sum_{1 \leq k \leq N_j}{(x_{k,j} e_{k,j}+\xi_{k,j} \eps_{k,j})},F\big{(}
\sum_{1 \leq k \leq N_j}{(x_{k,j} e_{k,j}+\xi_{k,j} \eps_{k,j})}\big)\Big).
\end{align*}
This implies that we can find symplectic coordinates
$$(x,\xi)=(x',x'';\xi',\xi'') \in \rr^{2n},$$
where $(x',\xi')$ and $(x'',\xi'')$ are some symplectic coordinates in $S^{\sigma \perp}$ and $S$
respectively, such that
\begin{equation}\label{inf11.5}\inc
q(x,\xi)=q_1(x',\xi')+q_2(x'',\xi''), \num
\end{equation}
with
\begin{equation}\label{inf12}\inc
q_1(x',\xi')=\sigma\big((x',\xi'),F|_{S^{\sigma \perp} }(x',\xi')\big) \num
\end{equation}
and
\begin{equation}\label{inf13}\inc
q_2(x'',\xi'')=  \sigma\big((x'',\xi''),F|_{S}(x'',\xi'')\big). \num
\end{equation}
Since from (\ref{h3}),
$$(\textrm{Re }F)S=\{0\},$$
the quadratic form $q_2$ is purely imaginary-valued and can be written as
\begin{equation}\label{sm14}\inc
q_2=i \tilde{q}_2, \num
\end{equation}
where $\tilde{q}_2$ is the real-valued quadratic form
\begin{equation}\label{sm15}\inc
\tilde{q}_2(x'',\xi'')= \sigma\big((x'',\xi''),\textrm{Im }F|_{S}(x'',\xi'')\big). \num
\end{equation}
When the additional assumption (\ref{la}) is fulfilled, this quadratic form $\tilde{q}_2$ must be elliptic on $\rr^{2n''}$. Since a real-valued elliptic quadratic form
is necessarily a positive definite or negative definite quadratic form, we deduce from the lemma 18.6.4 in \cite{hormander} that we can find new symplectic coordinates $(x'',\xi'')$ in $S$
and $\eps \in \{\pm 1\}$ such that
\begin{equation}\label{sm16}\inc
\tilde{q}_2(x'',\xi'')=\eps \sum_{j=1}^{n''}{\lambda_j(\xi_j''^2+x_j''^2)}, \num
\end{equation}
where $\lambda_j>0$ for all $j=1,...,n''$. This proves $(iii)$ in Proposition \ref{proposition1}.

Let us now study the properties of the quadratic form $q_1$. We denote by $F_1$ its Hamilton map
\begin{equation}\label{h8}\inc
F_1=F|_{S^{\sigma \perp}}, \num
\end{equation}
and define the following quadratic form
\begin{equation}\label{2.3.98}\inc
r(X')=\sum_{j=0}^{2n-1}{\textrm{Re } q_1\big((\textrm{Im } F_1)^j X'\big)}, \ X'=(x',\xi') \in S^{\sigma \perp}. \num
\end{equation}
Since from (\ref{sm4}), (\ref{inf11.5}) and (\ref{sm14}), $\textrm{Re } q_1$ is a non-positive quadratic form, we already know that $r$ is a non-positive quadratic form. We now prove that
$r$ is actually a negative definite quadratic form. Let us consider $X_0' \in S^{\sigma \perp}$ such that
$$r(X_0')=0.$$
The non-positivity of the quadratic form $\textrm{Re } q_1$ induces that for all $j=0,...,2n-1$,
\begin{equation}\label{2.3.99}\inc
\textrm{Re }q_1\big((\textrm{Im } F_1)^j X_0'\big)=0. \num
\end{equation}
Let us denote by $\textrm{Re }q_1(X';Y')$ the polarized form associated to $\textrm{Re }q_1$. We deduce from the Cauchy-Schwarz inequality, (\ref{10}) and (\ref{2.3.99}) that
for all $j=0,...,2n-1$ and $Y' \in S^{\sigma \perp}$,
\begin{align*}
|\textrm{Re }q_1\big(Y';(\textrm{Im }F_1)^j X_0'\big)|^2= & \ |\sigma\big(Y',\textrm{Re }F_1 (\textrm{Im }F_1)^j X_0'\big)|^2 \\
 \leq & \ [-\textrm{Re }q_1(Y')] [-\textrm{Re }q_1\big((\textrm{Im } F_1)^j X_0'\big) ]=0.
\end{align*}
It follows that for all $j=0,...,2n-1$ and $Y' \in S^{\sigma \perp}$,
$$\sigma\big(Y',\textrm{Re }F_1 (\textrm{Im }F_1)^j X_0'\big)=0,$$
which implies that for all $j=0,...,2n-1$,
\begin{equation}\label{2.3.100}\inc
\textrm{Re }F_1(\textrm{Im }F_1)^j X_0'=0, \num
\end{equation}
because from (\ref{2.3.89}), (\ref{h5}) and (\ref{h8}), $\textrm{Re }F_1(\textrm{Im }F_1)^j X_0' \in S^{\sigma \perp}$
and $S^{\sigma \perp}$ is a symplectic vector space.
Since $X_0' \in S^{\sigma \perp}$, we deduce from (\ref{h1bis}), (\ref{2.3.89}), (\ref{h5}), (\ref{h8}) and (\ref{2.3.100}) that $X_0' \in S \cap S^{\sigma \perp}=\{0\}$, which proves that $r$ is a negative definite quadratic form. This
proves $(ii)$ in Proposition \ref{proposition1}.

\bigskip

\noindent
\textit{Remark.} According to the previous proof, let us notice that the property $(ii)$ implies that for all $X \in \rr^{2n'}$, $X \neq 0$, there
exists $j_0 \in \{0,...,2n-1\}$ such that
\begin{equation}\label{sm17}\inc
\forall \ 0 \leq j \leq j_0-1, \ \textrm{Re }F_1(\textrm{Im }F_1)^j X=0, \ \textrm{Re }F_1(\textrm{Im }F_1)^{j_0} X \neq 0. \num
\end{equation}

\bigskip

Let us now prove that for all $T>0$, the average of the real part of the quadratic form $q_1$ by the flow defined by the Hamilton vector field of $\textrm{Im }q_1$,
$$\langle\textrm{Re }q_1\rangle_{T}(X')=\frac{1}{2T}\int_{-T}^{T}{\textrm{Re }q_1(e^{tH_{\textrm{Im}q_1}}X')dt},$$
is negative definite. Let us notice that this flow is globally defined since the symbol $\textrm{Im }q_1$ is quadratic.
Let us consider $X_0'$ in $\rr^{2n'}$ such that
\begin{equation}\label{sm18}\inc
\langle\textrm{Re }q_1\rangle_{T}(X_0')=0.\num
\end{equation}
Since $\textrm{Re }q_1$ is a non-positive quadratic form, it follows from (\ref{sm18}) that
$$\textrm{Re }q_1(e^{tH_{\textrm{Im}q_1}}X_0')=0,$$
for all $-T \leq t \leq T$. This implies in particular that for all $k \in \nn$,
\begin{equation}\label{h12}\inc
\frac{d^k}{dt^k}\big(\textrm{Re }q_1(e^{tH_{\textrm{Im}q_1}}X_0')\big)\big|_{t=0}=H_{\textrm{Im} q_1}^k\textrm{Re }q_1(X_0')=0.\num
\end{equation}
If $X_0' \neq 0$, we deduce from $(ii)$ that there exists $j_0 \in \{0,...,2n-1\}$ such that
\begin{equation}\label{st2}\inc
\forall \ 0 \leq j \leq j_0-1, \ \textrm{Re } q_1\big((\textrm{Im } F_1)^{j} X_0'\big)=0, \ \textrm{Re } q_1\big((\textrm{Im } F_1)^{j_0} X_0'\big)<0. \num
\end{equation}
Let us check that it would imply that
\begin{equation}\label{st4}\inc
H_{\textrm{Im} q_1}^{2j_0} \textrm{Re } q_1(X_0') \neq 0, \num
\end{equation}
and contradict (\ref{h12}). To prove (\ref{st4}), we use some arguments already used in \cite{mz} together with the following lemma also
proved in \cite{mz}.

\bigskip

\begin{lemma}\label{ll2}
If $q_1$ and $q_2$ are two complex-valued quadratic forms on $\rr^{2n}$, then the Hamilton map associated to the complex-valued quadratic form defined by the Poisson
bracket
$$\{q_1,q_2\}=\frac{\partial q_1}{\partial \xi}.\frac{\partial q_2}{\partial x}-\frac{\partial q_1}{\partial x}.\frac{\partial q_2}{\partial \xi},$$
is $-2[F_1,F_2]$ where $[F_1,F_2]$ stands for the commutator of $F_1$ and $F_2$, the Hamilton maps of $q_1$ and $q_2$.
\end{lemma}

\bigskip

We deduce from the previous lemma that the Hamilton map associated to the quadratic form $H_{\textrm{Im} q_1}^{2j_0} \textrm{Re } q_1$ is
\begin{equation}\label{st4.1}\inc
4^{j_0}[\textrm{Im }F_1,[\textrm{Im }F_1,[...,[\textrm{Im }F_1,\textrm{Re }F_1]...], \num
\end{equation}
with exactly $2j_0$ terms $\textrm{Im }F_1$ appearing in the formula. We can write
\inc
\begin{align*}\label{st4.5}
& \ 4^{j_0}[\textrm{Im }F_1,[\textrm{Im }F_1,[...,[\textrm{Im }F_1,\textrm{Re }F_1]...] \num \\
= & \ \sum_{j=0}^{2j_0}{(-1)^j c_j (\textrm{Im }F_1)^j \textrm{Re }F_1 (\textrm{Im }F_1)^{2j_0-j}},
\end{align*}
with $c_j >0$ for all $j=0,...,2j_0$. Indeed, by using the following identity
$$[P,[P,Q]]=P^2Q-2PQP+QP^2,$$
we can prove by induction that for all $n \in \nn^*$, there exist some positive constants $d_{n,j}$, $j=0,...,2n$ such that
$$[P,[P,[...,[P,Q]...]=\sum_{j=0}^{2n}{(-1)^j d_{n,j} P^j Q P^{2n-j}},$$
if there are exactly $2n$ terms $P$ in the left-hand side of the previous identity.
It follows from (\ref{st4.1}) and (\ref{st4.5}) that
\inc
\begin{align*}\label{st5}
& \ H_{\textrm{Im} q_1}^{2j_0} \textrm{Re } q_1(X_0')= (-1)^{j_0} c_{j_0}\sigma\big(X_0',(\textrm{Im }F_1)^{j_0}\textrm{Re }F_1(\textrm{Im }F_1)^{j_0}X_0' \big) \num \\
+ & \ \sum_{j=0}^{j_0-1}(-1)^j c_j \sigma\big(X_0',(\textrm{Im }F_1)^j\textrm{Re }F_1(\textrm{Im }F_1)^{2j_0-j}X_0' \big)\\
+ & \ \sum_{j=0}^{j_0-1}(-1)^{2j_0-j}c_{2j_0-j} \sigma\big(X_0',(\textrm{Im }F_1)^{2j_0-j}\textrm{Re }F_1(\textrm{Im }F_1)^{j}X_0' \big).
\end{align*}
Now on the one hand,
\begin{align*}
\sigma\big(X_0',(\textrm{Im }F_1)^{j_0}\textrm{Re }F_1(\textrm{Im }F_1)^{j_0}X_0' \big)= & \ (-1)^{j_0}\sigma\big((\textrm{Im }F_1)^{j_0}X_0',\textrm{Re }F_1(\textrm{Im }F_1)^{j_0}X_0' \big)\\
= & \ (-1)^{j_0}\textrm{Re }q_1\big((\textrm{Im }F_1)^{j_0}X_0'\big),
\end{align*}
by the skew-symmetry of the Hamilton map $\textrm{Im }F_1$. On the other hand, using (\ref{10}), (\ref{st2}),
and the Cauchy-Schwarz inequality, we get
\begin{align*}
& \ |\sigma\big(X_0',(\textrm{Im }F_1)^j\textrm{Re }F_1(\textrm{Im }F_1)^{2j_0-j}X_0' \big)|\\
=& \ |\sigma\big((\textrm{Im }F_1)^jX_0',\textrm{Re }F_1(\textrm{Im }F_1)^{2j_0-j}X_0' \big)|\\
= & \ -\textrm{Re }q_1\big((\textrm{Im }F_1)^jX_0';(\textrm{Im }F_1)^{2j_0-j}X_0'\big) \\
\leq & \ [-\textrm{Re }q_1((\textrm{Im }F_1)^jX_0')]^{\frac{1}{2}} [-\textrm{Re }q_1((\textrm{Im }F_1)^{2j_0-j}X_0')]^{\frac{1}{2}}=0
\end{align*}
and
\begin{align*}
& \ |\sigma\big(X_0',(\textrm{Im }F_1)^{2j_0-j}\textrm{Re }F_1(\textrm{Im }F_1)^{j}X_0' \big)|\\
=& \ |\sigma\big((\textrm{Im }F_1)^{2j_0-j}X_0',\textrm{Re }F_1(\textrm{Im }F_1)^{j}X_0' \big)|\\
= & \ -\textrm{Re }q_1\big((\textrm{Im }F_1)^{2j_0-j}X_0';(\textrm{Im }F_1)^{j}X_0'\big) \\
\leq & \ [-\textrm{Re }q_1((\textrm{Im }F_1)^{2j_0-j}X_0')]^{\frac{1}{2}} [-\textrm{Re }q_1((\textrm{Im }F_1)^{j}X_0')]^{\frac{1}{2}}=0
\end{align*}
if $j=0,...,j_0-1$. It follows from (\ref{st2}) and (\ref{st5}) that
$$H_{\textrm{Im} q_1}^{2j_0} \textrm{Re } q_1(X_0')=c_{j_0}\textrm{Re }q_1\big((\textrm{Im }F_1)^{j_0}X_0'\big) < 0,$$
because $c_{j_0} >0$. This proves (\ref{st4}) and ends the proof $(i)$. $\Box$

\bigskip

\noindent
\textit{Remark.} Let us notice that we have actually proved that the symbol~$q_1$ has a finite order $\tau$,
\begin{equation}\label{st1}\inc
1 \leq \tau \leq 4n-2, \num
\end{equation}
in every point of the set $q_1(\rr^{2n'}) \setminus \{0\}$.
We recall that the order $k(x_0,\xi_0)$ of a symbol $p(x,\xi)$ at a point $(x_0,\xi_0) \in \rr^{2n}$ (see section 27.2, chapter 27 in \cite{hormander}) is the element of
$\nn \cup \{+\infty\}$ defined by
\begin{equation}\label{t10.5}\inc
k(x_0,\xi_0)=\sup\big{\{}j \in \mathbb{Z} : p_I(x_0,\xi_0)=0, \ \forall\ 1 \leq |I| \leq j\big{\}}, \num
\end{equation}
where $I=(i_1,i_2,...,i_k) \in \{1,2\}^k$, $|I|=k$ and $p_I$ stands for the iterated Poisson brackets
$$p_{I}=H_{p_{i_1}}H_{p_{i_2}}...H_{p_{i_{k-1}}}p_{i_k},$$
where $p_1$ and $p_2$ are respectively the real and the imaginary part of the symbol $p$, $p=p_1+ip_2$. The order of a symbol $q$ at a point $z$ is then defined as the
maximal order of the symbol $p=q-z$ at every point $(x_0,\xi_0) \in \rr^{2n}$ verifying
$$p(x_0,\xi_0)=q(x_0,\xi_0)-z=0.$$

\section{Proofs of the main results}

\subsection{Heat semigroup smoothing for non-elliptic quadratic operators}\label{compact}
\init

In this section, we prove Theorem \ref{theorem1}. Let us consider a complex-valued quadratic form
$$q : \rr_x^n \times \rr_{\xi}^n \rightarrow \cc, \ n \in \nn^*,$$
with a non-positive real part
$$\textrm{Re }q \leq 0,$$
such that its singular space $S$ has a symplectic structure. We recall that this assumption is actually
fulfilled in the assumptions of Theorems \ref{theorem1}, \ref{theorem2} and
\ref{theorem3} since
this singular space is always symplectic when the symbol is elliptic on $S$.

We can then use the symplectic decomposition of the symbol obtained in section~\ref{sympl}.  We deduce from
(\ref{sm12}) and (\ref{sm13}) that there exists $\chi$, a real linear symplectic transformation of $\rr^{2n}$, such that
\begin{equation}\label{inv1}\inc
(q \circ \chi)(x,\xi)= q_1(x',\xi')+ iq_2(x'',\xi''), \ (x,\xi)=(x',x'';\xi',\xi'') \in \rr^{2n}, \num
\end{equation}
where $q_1$ is a complex-valued quadratic form on $\rr^{2n'}$ with a non-positive real part
\begin{equation}\label{inv2}\inc
\textrm{Re }q_1 \leq 0, \num
\end{equation}
and $q_2$ is a real-valued quadratic form verifying the properties stated in Proposition~\ref{proposition1}. The key point in our proof of Theorems \ref{theorem1}, \ref{theorem2}
and \ref{theorem3} is to
prove the following proposition.

\bigskip

\begin{proposition}\label{sperp}
If $n' \geq 1$, then the spectrum of the quadratic differential operator $q_1(x',\xi')^w$ is only composed of eigenvalues with finite multiplicity
$$
\sigma\big{(}q_1(x',\xi')^w\big{)}=\Big\{ \sum_{\substack{\lambda \in \sigma(F_1), \\  \emph{\textrm{Re}}(-i \lambda)<0} }
{\big{(}r_{\lambda}+2 k_{\lambda}
\big{)}(-i\lambda) : k_{\lambda} \in \nn}
\Big\},$$
where $F_1$ is the Hamilton map associated to the quadratic form $q_1$ and $r_{\lambda}$ is the
dimension of the space of generalized eigenvectors of $F_1$ in $\cc^{2n'}$
belonging to the eigenvalue $\lambda \in \cc$. Moreover, the
operator $q_1(x',\xi')^w$ generates a contraction semigroup such that
$$e^{tq_1(x',\xi')^w}u \in \mathcal{S}(\rr^{n'}),$$
for any $t>0$ and $u \in L^2(\rr^{n'})$.
\end{proposition}

\medskip
\noindent
{\it Remark.} It will be clear from the proof that Proposition 3.1.1 extends to the vector-valued case, so that if
${\mathcal H}$ is a complex Hilbert space and $u \in L^2(\rr^{n'};{\mathcal H})$ then for any $t>0$ we have
$e^{tq_1(x',\xi')^w}u \in \mathcal{S}(\rr^{n'};{\mathcal H})$.

\bigskip
\noindent

Theorem \ref{theorem1} directly follows from Proposition 3.1.1 together with the preceding remark. Indeed,
by using the symplectic invariance of the Weyl quantization given by the theorem~18.5.9 in
\cite{hormander},  we can find a metaplectic operator $U$, which is a unitary transformation on $\lde$ and an automorphism of $\mathcal{S}(\rr^n)$ such that
\begin{equation}\label{inv2.5}\inc
(q \circ \chi)(x,\xi)^w=U^{-1} q(x,\xi)^w U. \num
\end{equation}
This implies at the level of the generated semigroups that
\begin{equation}\label{inv3}\inc
e^{t(q \circ \chi)(x,\xi)^w}=U^{-1} e^{t q(x,\xi)^w} U, \ t \geq 0. \num
\end{equation}
Since from the tensorization of the variables (\ref{inv1}),
$$e^{t(q \circ \chi)(x,\xi)^w}=e^{tq_1(x',\xi')^w} e^{it q_2(x'',\xi'')^w},$$
we directly deduce from (\ref{inf12}), (\ref{inv3}), Proposition \ref{sperp} together with the following remark,
and the symplectic invariance of the Weyl quantization that if $(x',\xi')$ are some
symplectic coordinates on the symplectic space $S^{\sigma \perp}$ then for all $t >0$, $N \in \nn$ and
$u \in L^2(\rr^n)=L^2(\rr^{n'}; L^2(\rr^{n''}))$, we have
$$\big((1+|x'|^2+|\xi'|^2)^N\big)^we^{t q(x,\xi)^w}u \in L^2(\rr^n),$$
which proves Theorem \ref{theorem1}.

Let us now prove Proposition \ref{sperp}. For convenience, we drop the index and we consider a complex-valued
quadratic form
$$q : \rr_x^n \times \rr_{\xi}^n \rightarrow \cc, \ n \in \nn^*,$$
with a non-positive real part
\begin{equation}\label{kps2}\inc
\textrm{Re }q \leq 0, \num
\end{equation}
such that for all $T>0$, the average of the real part of the quadratic form $q$ by the flow defined by the Hamilton vector field of $\textrm{Im }q$,
\begin{equation}\label{inv5}\inc
\langle\textrm{Re }q\rangle_{T}(X)=\frac{1}{2T}\int_{-T}^{T}{\textrm{Re }q(e^{tH_{\textrm{Im}q}}X)dt}, \ X=(x,\xi) \in \rr^{2n}, \num
\end{equation}
is negative definite. We also know from (\ref{sm17}) that for all $X \in \rr^{2n}$, $X \neq 0$,
there exists $j_0 \in \nn$ verifying
\begin{equation}\label{inv6}\inc
\forall \ 0 \leq j \leq j_0-1, \ \textrm{Re }F(\textrm{Im }F)^j X=0, \ \textrm{Re }F(\textrm{Im }F)^{j_0} X \neq 0, \num
\end{equation}
if $F$ stands for the Hamilton map of the quadratic form $q$.

Let us denote by
\begin{equation}\inc
\label{eq3}
Q=q(x,\xi)^w,\num
\end{equation}
the quadratic differential operator defined by the Weyl quantization of the symbol~$q$.
When proving Proposition \ref{sperp}, we shall work with the metaplectic FBI-Bargmann
transform
\begin{equation}\inc
\label{eq4}
T u(x) = C \int_{\rr^n} e^{i\varphi(x,y)} u(y)\,dy,\ x\in \comp^n, \ C>0,\num
\end{equation}
where we may choose
$$
\varphi(x,y)=\frac{i}{2} (x-y)^2,
$$
as in the standard Bargmann transform.
Other quadratic phase functions $\varphi$ such that $\Im \varphi''_{y y}>0$ and $\det \varphi''_{x y}\neq 0$, are also
possible (see section 1 of~\cite{Sj95}).
It is well known that for a suitable choice of $C>0$, $T$ defines a unitary transformation
$$
T: L^2(\real^n)\rightarrow H_{\Phi_0}(\comp^n),
$$
where
\begin{equation}\inc
\label{eq4.1}
H_{\Phi_0}(\comp^n)={\rm Hol}(\comp^n)\cap L^2\big(\comp^n,
e^{-2\Phi_0(x)}L(dx)\big), \num
\end{equation}
with
$$
\Phi_0(x)=\sup_{y \in \rr^n} -\Im \varphi(x,y)=\frac{1}{2} \left(\Im x\right)^2,
$$
and $L(dx)$ being the Lebesgue measure in $\comp^n$.

\medskip
\noindent
{\it Remark.} Let us recall (see, e.g. section 3 of~\cite{SjZw}) that the same definitions apply in the vector-valued case,
so that we have a unitary operator
$$
T: L^2(\rr^n; {\mathcal H})\rightarrow H_{\Phi_0}(\comp^n; {\mathcal H}),
$$
where ${\mathcal H}$ is a complex Hilbert space.

\medskip
We recall next from~\cite{Sj95} that
\begin{equation}\inc
\label{eq5}
TQu=Q_0Tu,\ u\in \mathcal{S}(\rr^n), \num
\end{equation}
where $Q_0$ is a quadratic differential operator on $\comp^n$ whose
Weyl symbol $q_0$ satisfies
\begin{equation}\inc
\label{inv7}
q_0\circ \kappa_T = q. \num
\end{equation}
Here
\begin{equation}\inc
\label{eq5.1}
\kappa_T: \comp^{2n}\ni \big(y,-\varphi'_y(x,y)\big) \mapsto \big(x,\varphi'_x(x,y)\big)\in \comp^{2n}, \num
\end{equation}
is the complex linear canonical transformation associated to
$T$. From~\cite{Sj95}, we recall next that if we define
\begin{equation}\inc
\label{inv8}
\Lambda_{\Phi_0}=\Big\{\Big(x,\frac{2}{i}\frac{\partial \Phi_0}{\partial
x}(x)\Big): x\in \comp^n \Big\}, \num
\end{equation}
then we have
\begin{equation}\inc
\label{inv9}
\Lambda_{\Phi_0}=\kappa_T(\real^{2n}). \num
\end{equation}
When
$$\sigma = \sum_{j=1}^n d\xi_j \wedge dx_j,$$
is the complex symplectic (2,0)-form on $\comp^{2n}=\comp^n_x\times
\comp^n_{\xi}$, then the restriction $\sigma_{\Lambda_{\Phi_0}}$ of
$\sigma$ to $\Lambda_{\Phi_0}$ is real and nondegenerate. The map
$\kappa_T$ in (\ref{eq5.1}) can therefore be viewed as a canonical
transformation between the real symplectic spaces $\real^{2n}$ and
$\Lambda_{\Phi_0}$.

Continuing to follow~\cite{Sj95}, let us recall next that when realizing $Q_0$ as an unbounded operator on
$H_{\Phi_0}(\comp^n)$, we may use first the contour integral
representation
$$
Q_0 u(x) = \frac{1}{(2\pi)^n}\int_{\theta=\frac{2}{i}\frac{\partial \Phi_0}{\partial
x}\left(\frac{x+y}{2}\right)} e^{i(x-y)\cdot \theta}
q_0\Big(\frac{x+y}{2},\theta\Big)u(y)\,dy\,d\theta,
$$
and then, using that the symbol $q_0$ is holomorphic, by a contour
deformation we obtain the following formula for $Q_0$ as an unbounded operator on $H_{\Phi_0}(\comp^n)$,
\begin{equation}\inc
\label{eq5.5}
Q_0 u(x) = \frac{1}{(2\pi)^n}\int_{\theta=\frac{2}{i}\frac{\partial \Phi_0}{\partial
x}\left(\frac{x+y}{2}\right)+it\overline{(x-y)}} e^{i(x-y)\cdot \theta}
q_0\Big(\frac{x+y}{2},\theta\Big)u(y)\,dy\,d\theta, \num
\end{equation}
for any $t>0$.

We shall now discuss certain IR-deformations of the real phase space
$\real^{2n}$, where the averaging procedure along the flow defined by the Hamilton vector field of $\textrm{Im }q$ (see (\ref{inv5}))
plays an important r\^ole. To that end, let $G=G_T$ be a real-valued quadratic form on $\real^{2n}$ such that
\begin{equation}\inc
\label{eq6}
H_{\textrm{Im} q} G = -\Re q+\langle{\Re q}\rangle_T. \num
\end{equation}
As in~\cite{HeHiSj}, we solve (\ref{eq6}) by setting
\begin{equation}\inc
\label{eq7}
G(X)=\int_{\rr} k_T(t) \Re q(e^{tH_{\textrm{Im} q}}X)dt,\num
\end{equation}
where $k_T(t)=k(t/2T)$ and $k\in C(\real\backslash\{0\})$ is the odd function given by
$$
k(t)=0 \textrm{ for } |t| \geq \frac{1}{2} \textrm{ and } k'(t)=-1 \textrm{ for } 0<|t| < \frac{1}{2}.
$$
Let us notice that $k$ and $k_T$ have a jump of size $1$ at the origin.
Associated with $G$ there is a linear IR-manifold, defined for $0\leq \eps\leq \eps_0$, with $\eps_0>0$
small enough,
\begin{equation}\inc
\label{eq8}
\Lambda_{\eps G}=e^{i\eps H_G}(\rr^{2n}) \subset \comp^{2n},\num
\end{equation}
where $e^{i\eps H_G}$ stands for the flow generated by the linear Hamilton vector field $i\eps H_G$ taken at the time 1.
It is then well-known and easily checked (see, for instance, sections 3 and 5 in~\cite{HeSjSt}), that
\begin{equation}\inc
\label{eq9}
\kappa_T(\Lambda_{\eps
G})=\Lambda_{\tilde{\Phi}_{\eps}}:=\Big\{\Big(x,\frac{2}{i}\frac{\partial
\tilde{\Phi}_{\eps}}{\partial x}(x)\Big) : x\in \comp^n\Big\},\num
\end{equation}
where $\tilde{\Phi}_{\eps}$ is a strictly plurisubharmonic quadratic form on $\comp^n$, such that
\begin{equation}\inc\label{lin1}
\tilde{\Phi}_{\eps}(x)=\Phi_0(x)+\eps G(\Re x, -\Im x)+\mathcal{O}(\eps^2\abs{x}^2). \num
\end{equation}
Associated with the function $\tilde{\Phi}_{\eps}$ is the weighted space of
holomorphic functions $H_{\tilde{\Phi}_{\eps}}(\comp^n)$
defined as in
(\ref{eq4.1}). The operator $Q_0$ can also be defined as an unbounded operator
$$
Q_0: H_{\tilde{\Phi}_{\eps}}(\comp^n) \rightarrow H_{\tilde{\Phi}_{\eps}}(\comp^n),
$$
if we make a new contour deformation in (\ref{eq5.5}) and set
\begin{equation}\inc
\label{eq10}
Q_0 u(x) = \frac{1}{(2\pi)^n}\int_{\theta=\frac{2}{i}\frac{\partial \tilde{\Phi}_{\eps}}{\partial
x}\left(\frac{x+y}{2}\right)+it\overline{(x-y)}} e^{i(x-y)\cdot \theta}
q_0\Big(\frac{x+y}{2},\theta\Big)u(y)\,dy\,d\theta,\num
\end{equation}
for any $t>0$. By coming back to the real side by the FBI-Bargmann transform, the operator $Q_0$ can be viewed as an
unbounded operator on $L^2(\rr^n)$ with the Weyl symbol
\begin{equation}\inc\label{esk1}
\widetilde{q}(X)=q\big(e^{i \eps H_G}X\big), \num
\end{equation}
and here the real part of this expression is easily seen to be equal to
$$
\textrm{Re }\widetilde{q}(X)=\textrm{Re }q(X) +\eps H_{\textrm{Im}q} G(X)+\mathcal{O}(\eps^2\abs{X}^2).
$$
It follows therefore from (\ref{kps2}), (\ref{eq6}) and the assumption that the quadratic form (\ref{inv5}) is negative definite,
that
\begin{equation}\inc
\label{eq10.1}
-\textrm{Re } \widetilde{q}(X)\geq \frac{\eps}{C}\abs{X}^2, \ C>1,\ X \in \real^{2n}, \num
\end{equation}
for $0<\eps \ll 1$. We may therefore apply Theorem 3.5 of~\cite{sjostrand} to the
operator $Q_0$ viewed as an unbounded operator on $H_{\tilde{\Phi}_{\eps}}(\cc^n)$.

\bigskip

\begin{lemma}\label{spec}
Let us consider $Q_0$ as an unbounded operator on $H_{\tilde{\Phi}_{\eps}}(\cc^n)$, for $0<\eps \leq \eps_0$, with $\eps_0>0$
sufficiently small. The spectrum of the operator $Q_0$ is only composed of eigenvalues with finite multiplicity
\begin{equation}\inc
\label{eq11}
\sigma(Q_0)=\Big\{ \sum_{\substack{\lambda\in \sigma(F), \\ \emph{\textrm{Re}}(-i\lambda)<0}}
\left(r_{\lambda}+2 k_{\lambda}\right)(-i\lambda) : \ k_{\lambda}\in
\nat \Big\},\num
\end{equation}
where $F$ is the Hamilton map associated to the quadratic form $q$ and $r_{\lambda}$ is the dimension of the space of generalized
eigenvectors of $F$ in $\cc^{2n}$ belonging to the eigenvalue $\lambda\in \comp$.
\end{lemma}

\bigskip

To get the statement of this lemma, it suffices to combine Theorem 3.5 of~\cite{sjostrand} together with the observation
that the Hamilton maps $F$ and $\widetilde{F}$ of the quadratic forms $q$ and $\widetilde{q}$, respectively, are isospectral since from (\ref{esk1}) the symbols $q$ and
$\widetilde{q}$ are related by a canonical transformation.

Having determined the spectrum of $Q_0$ in the weighted space $H_{\tilde{\Phi}_{\eps}}(\cc^n)$, $0<\eps\ll 1$,
we now come to the proof of Proposition~\ref{sperp}.
In doing so, we shall study the spectral properties of the holomorphic quadratic differential operator $Q_0$ acting on $H_{\Phi_0}(\cc^n)$.

We shall consider the heat evolution equation associated to the operator $Q_0$. Let us notice explicitly that we got this idea by studying Remark 11.7
in~\cite{HeSjSt}, and indeed, the following argument can be seen as a natural continuation of some ideas sketched in that remark. Using
Fourier integral operators with quadratic phase in the complex domain, we may describe the heat semigroup $e^{tQ_0}$ for $0\leq t \leq t_0$,
when $t_0>0$ is small enough. More precisely, we are interested in
solving
$$
\left\lbrace
\begin{array}{c}
\displaystyle   \frac{\partial u}{\partial t}(t,x)-Q_0 u(t,x)=0 \\
u(t,\textrm{\textperiodcentered})|_{t=0}=u_0 \in H_{\Phi_0}(\cc^n).
\end{array} \right.
$$
Let $\varphi(t,x,\eta)$ be the quadratic form in the variables $x$, $\eta$, depending
smoothly on $t$, $0\leq t\leq t_0\ll 1$, and solving the Hamilton-Jacobi equation
$$
\left\lbrace
\begin{array}{c}
\displaystyle i\frac{\partial \varphi}{\partial t}(t,x,\eta)-q_0\Big(x,\frac{\partial \varphi}{\partial x}(t,x,\eta)\Big)=0 \\
\varphi(t,x,\eta)|_{t=0}=x\cdot \eta.
\end{array} \right.
$$
We know that for $0\leq t \leq t_0 \ll 1$, $\varphi(t,x,\eta)$ can be obtained as a generating
function of the complex canonical transformation
$$e^{itH_{q_0}} : \big(\varphi'_{\eta}(t,x,\eta),\eta\big)\mapsto \big(x,\varphi'_x(t,x,\eta)\big).$$
Then for $t\geq 0$ small enough, the operator $e^{tQ_0}$ acting on $H_{\Phi_0}(\cc^n)$ has the form
$$e^{tQ_0} u = \frac{1}{(2\pi)^n}\int_{\Gamma_x} e^{i(\varphi(t,x,\eta)-y\cdot \eta)} a(t,x,y,\eta) u(y) dy\, d\eta,$$
where $a(t,x,y,\eta)$ is a suitable amplitude which we need not specify here, and, following the
general theory of~\cite{Sj82}, we take $\Gamma_x$ to be a
suitable contour passing through the critical point of the function
$$(y,\eta)\mapsto -\textrm{Im}\big(\varphi(t,x,\eta)-y\cdot\eta\big)+\Phi_0(y).$$
We then know from the general theory that the operator $e^{tQ_0}$ is bounded
\begin{equation}\inc
\label{invb1}
e^{tQ_0}: H_{\Phi_0}(\cc^n) \rightarrow H_{\Phi_t}(\cc^n),\num
\end{equation}
where $\Phi_t$ is a strictly plurisubharmonic quadratic form on
$\comp^n$, depending smoothly on $t$, such that if
\begin{equation}\inc\label{kps40}
\Lambda_{\Phi_t}=\Big\{\Big(x,\frac{2}{i}\frac{\partial
\Phi_t}{\partial x}(x)\Big) : x\in \comp^n \Big\}, \num
\end{equation}
then
\begin{equation}\inc
\label{eq20}
\Lambda_{\Phi_t} =
\exp\big(t\widehat{H_{-\frac{1}{i}q_0}}\big)\left(\Lambda_{\Phi_0}\right). \num
\end{equation}
Here, when $f$ is a holomorphic function on $\comp^{2n}=\comp^n_x\times \comp^n_{\xi}$, $H_f$ is the
standard Hamilton field of $f$, of type (1,0), given by the usual formula
$$
H_f = \sum_{j=1}^n \Big(\frac{\partial f}{\partial
\xi_j}\frac{\partial}{\partial x_j}-\frac{\partial f}{\partial
x_j}\frac{\partial}{\partial \xi_j}\Big),
$$
and
$\widehat{H_f}=H_f+\overline{H_f}$ is the corresponding real vector
field.

It follows from the classical Hamilton-Jacobi theory applied with respect to the real symplectic form $\Im \sigma$, where
$$\sigma=\sum_{j=1}^n d\xi_j \wedge dx_j,$$
is the complex symplectic (2,0)-form on $\comp^{2n}$,
that the quadratic form $\Phi(t,x)=\Phi_t(x)$ introduced in
(\ref{invb1}), satisfies the eikonal equation
\begin{equation}\inc
\label{eikonal}
\left\lbrace
\begin{array}{c}
\displaystyle \frac{\partial \Phi}{\partial t}(t,x)- \Re \Big[q_0\Big(x,\frac{2}{i}\frac{\partial \Phi}{\partial x}(t,x)\Big)\Big]=0 \\
\Phi(t,\textrm{\textperiodcentered})|_{t=0}=\Phi_0.
\end{array} \right.
\num
\end{equation}
See also~\cite{SjHokk} for a detailed (and much more general)
discussion of this point.

\medskip
Instrumental in the proof of Proposition~\ref{sperp} is the following result.

\bigskip

\begin{lemma}\label{phase}
For each $T_0>0$
small enough, there exists $\alpha=\alpha(T_0)>0$ such that
\begin{equation}\inc
\label{eq20.1}
\Phi_{T_0}(x)\leq \Phi_0(x) - \alpha \abs{x}^2,\ x\in \comp^n. \num
\end{equation}
\end{lemma}

\bigskip

Once Lemma \ref{phase} has been established, it is easy to finish the proof of the first result in Proposition~\ref{sperp}. Indeed, elementary arguments together with
(\ref{eq4.1}) and (\ref{eq20.1}) show that the natural embedding
$H_{\Phi_t}(\cc^n) \rightarrow H_{\Phi_0}(\cc^n)$
is compact for $t>0$ small, and hence by using the semigroup property, we deduce that the semigroup
\begin{equation}\inc\label{kps41}
e^{tQ_0}: H_{\Phi_0}(\cc^n) \rightarrow H_{\Phi_0}(\cc^n),\num
\end{equation}
is compact for each $t>0$. An application of Theorem 2.20 in~\cite{davies} then shows that the spectrum of the operator $Q_0$ acting on $H_{\Phi_0}(\cc^n)$ consists of a countable
discrete set of eigenvalues each of finite multiplicity.

When deriving the explicit description of the spectrum of $Q_0$ on $H_{\Phi_0}(\cc^n)$, we argue in the following way. Let us
assume that $\lambda \in \comp$ is an eigenvalue of $Q_0$ on $H_{\Phi_0}(\cc^n)$
and let $u_0 \in H_{\Phi_0}(\cc^n)$ be a corresponding eigenvector,
$$Q_0 u_0 =\lambda u_0.$$
We deduce from (\ref{invb1}) that
$$e^{tQ_0} u_0 \in H_{\Phi_t}(\cc^n),$$
and since
$$e^{tQ_0} u_0 = e^{t\lambda} u_0,$$
it follows from Lemma \ref{phase} that
$$u_0 \in H_{\Phi_0-\delta\abs{x}^2}(\cc^n),$$
for some $\delta>0$. In particular, we obtain from (\ref{lin1}) that $u_0\in H_{\tilde{\Phi}_{\eps}}(\cc^n)$ for $\eps>0$ small enough, and hence that $\lambda$ is in the spectrum of
the operator $Q_0$
acting on $H_{\tilde{\Phi}_{\eps}}(\cc^n)$, which has been described in Lemma \ref{spec}.

On the other hand, if $\lambda$ is in the spectrum of $Q_0$ acting
on $H_{\tilde{\Phi}_{\eps}}(\cc^n)$ and $u_0 \in H_{\tilde{\Phi}_{\eps}}(\cc^n)$ is a corresponding eigenvector, with $\eps>0$ sufficiently small, then we have
$$e^{tQ_0}u_0=e^{t\lambda}u_0\in H_{\widetilde{\Phi}_{\eps,t}}(\cc^n),$$
where $\widetilde{\Phi}_{\eps,t}$ is a quadratic form on $\comp^n$ depending smoothly on $t\geq 0$ and $\eps\geq 0$, for $\eps$
sufficiently small, which satisfies the eikonal equation (\ref{eikonal}) along with the initial condition
$$\widetilde{\Phi}_{\eps,t}(x)|_{t=0}=\widetilde{\Phi}_{\eps}(x).$$
It follows from (\ref{lin1}) and (\ref{eikonal}) that
$$\widetilde{\Phi}_{\eps,t}(x)=\Phi_t(x)+\mathcal{O}(\eps \abs{x}^2),$$
where the implicit constant is uniform in $0\leq t\leq t_0$, for $t_0>0$ small enough.
By taking $T_0>0$ small enough but fixed such that $0<T_0 \leq t_0$ and (\ref{eq20.1}) holds, we get
$$
u_0= e^{-T_0\lambda} e^{T_0 Q} u_0 \in H_{\tilde{\Phi}_{\eps,T_0}}(\cc^n)=
H_{\Phi_{T_0}+\mathcal{O}(\eps\abs{x}^2)}(\cc^n).
$$
In view of (\ref{eq20.1}), we can choose $\eps_0>0$ small enough such that for all $0<\eps \leq \eps_0$,
$$u_0 \in H_{\Phi_0-\tilde{\delta}|x|^2}(\cc^n) \subset H_{\Phi_0}(\cc^n),$$
where $\tilde{\delta}$ is a positive constant.
It follows that $\lambda$ is also in the spectrum of the operator
$Q_0$ acting on $H_{\Phi_0}(\cc^n)$. Altogether, this shows that the spectrum of $Q_0$ acting on $H_{\Phi_0}(\cc^n)$ is equal to
the spectrum of $Q_0$ acting on $H_{\widetilde{\Phi}_{\eps}}(\cc^n)$, for $\eps>0$ sufficiently small, and furthermore, that the
algebraic multiplicities agree. We have therefore identified the spectrum of $Q_0$ on $H_{\Phi_0}(\cc^n)$ and also the spectrum
of the operator $Q$ on $L^2(\rr^n)$ by coming back to the real side.

\bigskip

\noindent {\it Remark}. In the argument above we have worked with the eigenfunctions of the operator $Q_0$ acting on
$H_{\widetilde{\Phi}_{\eps}}(\cc^n)$ for some sufficiently small but fixed $\eps>0$. It may be interesting to notice that the
(generalized) eigenfunctions of the operator $Q_0$ do not depend on $\eps>0$, for $0<\eps\leq \eps_0$, with $\eps_0>0$ small enough.
See also Remark~11.7 in~\cite{HeSjSt}. While refraining from providing a
detailed proof of this statement, let us mention that its validity
relies crucially upon the fact that the quadratic symbol $q_0$ is
elliptic on $\Lambda_{\widetilde{\Phi}_{\eps}}$, $0<\eps\leq
\eps_0$. In the terminology of~\cite{Sj82}, this is an instance of
the principle of non-characteristic deformations.

\bigskip

We shall now prove Lemma \ref{phase}. Integrating (\ref{eikonal}) for $t=0$ to $t=T$, we obtain
$$\Phi_T(x)-\Phi_0(x)=\int_0^{T} \Re
\Big[q_0\Big(x,\frac{2}{i}\frac{\partial \Phi}{\partial
x}(t,x)\Big)\Big]\,dt,\ T>0.$$
Here the integral in the right hand side is a real-valued quadratic form on $\comp^n$,
and Lemma~\ref{phase} is implied by the claim that it is negative definite. When proving the claim, we shall use the following general
relation, explained for instance in~\cite{MelinSj},
\begin{equation}\inc
\label{eq21}
\widehat{H_{if}} = H^{-{\rm Im}\sigma}_{-\textrm{Re}f},\num
\end{equation}
valid for a holomorphic function $f$ on $\comp^{2n}$.

It follows from (\ref{eq21}) that $\Re q_0$ is constant along the flow of the Hamilton vector field
$$\widehat{H_{-\frac{1}{i}q_0}}=H^{-\textrm{Im}\sigma}_{-\textrm{Re }q_0}.$$
We therefore deduce from (\ref{kps2}), (\ref{inv7}), (\ref{inv9}) and (\ref{eq20}) that
\begin{equation}\inc \label{kps30}
\textrm{Re}\big[q_0|_{\Lambda_{\Phi_t}}\big] \leq 0. \num
\end{equation}
According to (\ref{kps40}), the equation (\ref{eikonal}) implies therefore that
$$\frac{\partial \Phi}{\partial t}(t,x)\leq 0,$$
so that the function
$$t\mapsto \Phi_t(x),$$
is decreasing. Let us assume that there exists $x_0\in \comp^n$, $x_0\neq 0$,
such that
$$
\int_0^T \Re
\Big[q_0\Big(x_0,\frac{2}{i}\frac{\partial \Phi}{\partial
x}(t,x_0)\Big)\Big]\,dt=0.
$$
Since from (\ref{kps40}) and (\ref{kps30}), the integrand is non-positive, it follows that for all $0\leq
t\leq T$,
\begin{equation}\inc \label{kps31}
\Re \Big[ q_0\Big(x_0,\frac{2}{i}\frac{\partial \Phi}{\partial
x}(t,x_0)\Big)\Big]=0. \num
\end{equation}
Therefore, in view of (\ref{eikonal}), we get
$$\frac{\partial \Phi}{\partial t}(t,x_0)=0,$$
for all $0\leq t \leq T$. Here, the quadratic form
$$f_t(x)=\frac{\partial \Phi}{\partial t}(t,x)=\Re \Big[ q_0\Big(x,\frac{2}{i}\frac{\partial \Phi}{\partial x}(t,x)\Big)\Big] \leq 0,$$
is non-positive and such that $f_t(x_0)=0$ for all $0 \leq t \leq T$. It follows that
$$\nabla_{{\rm Re} x, {\rm Im} x}f_t(x_0)=0,$$
for all $0\leq t\leq T,$ and therefore
$$
\frac{\partial f_t}{\partial x}(x_0)=\frac{\partial^2 \Phi}{\partial x
\partial t}(t,x_0)=0,
$$
for all $0\leq t\leq T$.
Hence the function
$$
t\mapsto \frac{\partial \Phi}{\partial x}(t,x_0),
$$
does not depend on $t$ for $0\leq t\leq T$, so that
\begin{equation}\inc \label{kps32}
\frac{\partial \Phi}{\partial x}(t,x_0)=\frac{\partial \Phi}{\partial
x}(0,x_0)=\frac{\partial \Phi_0}{\partial x}(x_0), \num
\end{equation}
for all $0\leq t\leq T$.
Since from (\ref{kps40}), the point
\begin{equation}\inc \label{kps33}
\Big(x_0,\frac{2}{i}\frac{\partial \Phi}{\partial
x}(t,x_0)\Big)=\Big(x_0,\frac{2}{i}\frac{\partial \Phi_0}{\partial
x}(x_0)\Big), \num
\end{equation}
belongs to $\Lambda_{\Phi_t}$ for all $0\leq t\leq T$, we obtain from (\ref{inv8}) and (\ref{eq20}) that there exists $y_0(t)\in \comp^n$ such that
\begin{equation}\inc \label{kps34}
\Big(x_0,\frac{2}{i}\frac{\partial \Phi}{\partial
x}(t,x_0)\Big)=\exp(t\widehat{H_{-\frac{1}{i}q_0}})\Big(y_0(t), \frac{2}{i}\frac{\partial \Phi_0}{\partial
x}\big(y_0(t)\big)\Big). \num
\end{equation}
It follows from (\ref{kps31}) that for all $0\leq t\leq T$,
\begin{equation}\inc\label{kps35}
\textrm{Re}\Big[q_0\Big(x_0,\frac{2}{i}\frac{\partial
\Phi}{\partial x}(t,x_0)\Big)\Big]=\textrm{Re}\Big[q_0\Big(y_0(t),\frac{2}{i}\frac{\partial \Phi_0}{\partial
x}\big(y_0(t)\big)\Big)\Big]=0, \num
\end{equation}
because $\textrm{Re }q_0$ is constant along the flow of the Hamilton vector field
$$\widehat{H_{-\frac{1}{i}q_0}} = H^{-{\rm Im}\sigma}_{-\textrm{Re}q_0}.$$
Let us now set
\begin{equation}\label{kps39}\inc
L_0=\big\{\tilde{X} \in \Lambda_{\Phi_0} : \textrm{Re}[q_0(\tilde{X})]=0 \big\}.\num
\end{equation}
We can notice by using similar arguments as in (\ref{2.3.99}) and (\ref{2.3.100}) that
$$\{X \in \rr^{2n}, \textrm{Re }\tilde{q}(X)=0\}=\textrm{Ker}(\textrm{Re }\tilde{F}) \cap \rr^{2n},$$
for any complex-valued quadratic form $\tilde{q}$ with a non-positive real part if $\textrm{Re }\tilde{F}$ is the Hamilton map of the quadratic form $\textrm{Re }\tilde{q}$.
We therefore deduce from (\ref{inv7}) and (\ref{inv9}) that
\begin{equation}\inc
\label{kps40yyy}
L_0=\kappa_T\big({\rm Ker}(\Re F)\cap \real^{2n}\big). \num
\end{equation}
We get from (\ref{inv8}), (\ref{kps33}), (\ref{kps34}), (\ref{kps35}), (\ref{kps39}) and (\ref{kps40yyy}) that
\begin{equation}\label{kps35bis}\inc
\Big(y_0(t),\frac{2}{i}\frac{\partial \Phi_0}{\partial
x}\big(y_0(t)\big)\Big)=\exp(t\widehat{H_{\frac{1}{i}q_0}})\Big(x_0,\frac{2}{i}\frac{\partial \Phi_0}{\partial
x}(x_0)\Big) \in L_0, \num
\end{equation}
for all $0\leq t\leq T$,
and therefore
\begin{equation}\inc
\label{eq25}
\Re F \Big(\kappa_T^{-1}\Big(\exp(t\widehat{H_{\frac{1}{i}q_0}})\Big(x_0,\frac{2}{i}\frac{\partial
\Phi_0}{\partial x}(x_0)\Big)\Big)\Big)=0, \num
\end{equation}
for all $0\leq t\leq T$. In view of (\ref{inv8}) and (\ref{inv9}), we may write
\begin{equation}\inc
\label{kps36}
\Big(x_0,\frac{2}{i}\frac{\partial \Phi_0}{\partial
x}(x_0)\Big)=\kappa_T(X_0), \num
\end{equation}
for some $X_0\in \real^{2n}$, $X_0\neq 0$ because $x_0 \neq 0$. We can now deduce from (\ref{inv6}) that there exists an integer $m \in \nn$ such that
\begin{equation}\inc
\label{eq26}
\Re F\left(\Im F\right)^j X_0=0,\quad 0\leq j<m, \num
\end{equation}
while
\begin{equation}\inc
\label{eq27}
\Re F\left(\Im F\right)^m X_0\neq 0.\num
\end{equation}
On the other hand, we get from (\ref{eq25}) and (\ref{kps36}) that
\begin{equation}\inc
\label{eq27.5}
\Re F\Big(\kappa_T^{-1}
\big(\widehat{H_{\frac{1}{i}q_0}}\big)^j\kappa_T(X_0)\Big)=0, \num
\end{equation}
for all $j \in \nn$ because $q_0$ is a quadratic form. We shall establish the following result.

\bigskip

\begin{lemma}\label{lem3.1.3}
For all $0\leq j\leq m$, we have
$$
\big(\widehat{H_{\frac{1}{i}q_0}}\big)^j
\kappa_T(X_0)=\big(H^{\sigma_{\Lambda_{\Phi_0}}}_{{\rm
Im}\,q_0}\big)^j \kappa_T(X_0).
$$
\end{lemma}

\bigskip

\begin{proof}
When proving this lemma, we shall argue by induction with respect to
$j$, and start with the case $j=0$, which is of course fulfilled. Let us
recall from (\ref{kps35bis}) and (\ref{kps36}) that $\kappa_T(X_0)\in L_0$, and notice that, as recalled
for example in section 11 (Remark 11.7) in \cite{HeSjSt} (see also~\cite{MelinSj}), we
have at the points of $L_0$,
\begin{equation}\inc
\label{eqSj}
\widehat{H_{\frac{1}{i}q_0}}=H^{\sigma_{\Lambda_{\Phi_0}}}_{{\rm
Im}\,q_0}. \num
\end{equation}
Let us now check that for all $0\leq j\leq m-1$,
\begin{equation}\inc
\label{eq28}
\big(H^{\sigma_{\Lambda_{\Phi_0}}}_{{\rm Im}\,q_0}\big)^j
\kappa_T(X_0)\in L_0.\num
\end{equation}
Let us consider $0 \leq j \leq m-1$ and notice from (\ref{inv7}) that
\begin{equation}\label{kps37}\inc
\kappa_T(H_{\textrm{Im }q}) \kappa_T^{-1}= H_{\textrm{Im }q_0}^{\sigma_{\Lambda_{\Phi_0}}}.\num
\end{equation}
Since a direct computation using (\ref{10}) shows that
\begin{equation}\label{kps38}\inc
H_{\textrm{Im }q} =2 \textrm{Im }F, \num
\end{equation}
we obtain by using (\ref{10}), (\ref{inv7}), (\ref{eq26}) and (\ref{kps37}) that
\begin{align*}
& \ \textrm{Re}\Big[q_0\big((H^{\sigma_{\Lambda_{\Phi_0}}}_{{\rm
Im}\,q_0})^j \kappa_T(X_0)\big)\Big]=\textrm{Re}\Big[q_0\big(\kappa_T H_{\Im q}^j (X_0)\big)\Big]\\
=& \ \Re q\left(2^j (\Im F)^j X_0\right)
=2^{2j}\sigma\big((\textrm{Im }F)^jX_0, \textrm{Re }F(\textrm{Im }F)^jX_0 \big)=0,
\end{align*}
for any $0 \leq j \leq m-1$.
Thus we have
verified (\ref{eq28}), and by an application of (\ref{eqSj}), we get that for all $0\leq j\leq m-1$,
$$
\widehat{H_{\frac{1}{i}q_0}}\big((\widehat{H_{\frac{1}{i}q_0}})^j
\kappa_T(X_0)\big)=H^{\sigma_{\Lambda_{\Phi_0}}}_{{\rm
Im}\,q_0}\big((\widehat{H_{\frac{1}{i}q_0}})^j
\kappa_T(X_0)\big)=(H^{\sigma_{\Lambda_{\Phi_0}}}_{{\rm Im}\,
q_0})^{j+1}\kappa_T(X_0),
$$
if
$$(\widehat{H_{\frac{1}{i}q_0}})^j\kappa_T(X_0)=(H^{\sigma_{\Lambda_{\Phi_0}}}_{{\rm Im}\,
q_0})^{j}\kappa_T(X_0).$$
This proves by induction Lemma~\ref{lem3.1.3}.
\end{proof}

It is now easy to finish the proof of Lemma~\ref{phase}. By using
(\ref{eq27.5}) when $j=m$, (\ref{kps37}), (\ref{kps38}) and applying Lemma~\ref{lem3.1.3}, we get
$$
0=\Re F\big(\kappa_T^{-1} (H^{\sigma_{\Lambda_{\Phi_0}}}_{{\rm
Im}\,q_0})^m \kappa_T(X_0)\big)=\Re F ( H_{{\rm Im}\,q}^m
X_0)=2^m \Re F(\Im F)^m X_0,
$$
which contradicts (\ref{eq27}) and completes the proof of Lemma~\ref{phase}.

Let us finally notice that the semigroup $e^{tQ}$, $t>0$, is strongly regularizing on $L^2(\rr^n)$. We actually deduce from (\ref{invb1}), (\ref{eq20.1}) and the fundamental property of semigroups
that for all $t>0$, there exists $\delta >0$ such that
$$ \forall u \in H_{\Phi_0}(\cc^n), \ e^{tQ_0} u \in H_{\Phi_0-\delta\abs{x}^2}(\cc^n),$$
on the FBI transform side. By using the fact that a holomorphic function $U$ on
$\comp^n$ is of the form $Tu$ for some $u\in \mathcal{S}(\rr^n)$, if and only if
$$
\forall N \in \nn, \ \int_{\cc^n} \abs{U(x)}^2 e^{-2\Phi_0(x)} \langle{x}\rangle^N\, L(dx)<+\infty,
$$
(see for instance \cite{Sj95}) we finally obtain that
$$\forall t>0, \forall u \in L^2(\rr^n), \ e^{tQ}u \in \mathcal{S}(\rr^n),$$
which ends the proof of Proposition~\ref{sperp}.

\subsection{Large time behavior of contraction semigroups}\label{beha}
\init

In this section, we prove Theorem \ref{theorem3}. Let us consider a complex-valued quadratic form
$$q : \rr_x^n \times \rr_{\xi}^n \rightarrow \cc, \ n \in \nn^*,$$
with a non-positive real part
$$\textrm{Re }q \leq 0,$$
such that its singular space $S$ has a symplectic structure.

Let us assume that the real part of the symbol $q$ is a non-zero quadratic form
$$\exists X_0 \in \rr^{2n}, \ \textrm{Re }q(X_0) \neq 0.$$
This implies that the singular space is distinct from the whole phase space $S \neq \rr^{2n}$ because from (\ref{10}) and (\ref{h1}),
$$\forall X \in S, \ \textrm{Re }q(X)=\sigma(X,\textrm{Re }FX)=0.$$
It proves that $(ii)$ implies $(iii)$ in Theorem \ref{theorem3}.

Let us now assume that $S \neq \rr^{2n}$ and prove $(i)$. We deduce from (\ref{sm12}) and (\ref{inf12}) that there exists
$\chi$, a real linear symplectic transformation of $\rr^{2n}$, such that
\begin{equation}\label{decay1}\inc
(q \circ \chi)(x,\xi)= q_1(x',\xi')+ iq_2(x'',\xi''), \ (x,\xi)=(x',x'';\xi',\xi'') \in \rr^{2n}, \num
\end{equation}
where $q_1$ is a complex-valued quadratic form on $\rr^{2n'}$, $n' \geq 1$, with a non-positive real part
and $q_2$ is a real-valued quadratic form verifying the properties stated in Proposition \ref{proposition1}.
By using the symplectic invariance of the Weyl quantization given by the theorem~18.5.9 in
\cite{hormander},  we can find a metaplectic operator $U$, which is a unitary transformation on $\lde$ and an automorphism of $\mathcal{S}(\rr^n)$ such that
$$(q \circ \chi)(x,\xi)^w=U^{-1} q(x,\xi)^w U.$$
This implies at the level of the generated semigroups that
$$e^{t(q \circ \chi)(x,\xi)^w}=U^{-1} e^{t q(x,\xi)^w} U, \ t \geq 0$$
and
\begin{equation}\label{decay2}\inc
\|e^{t(q \circ \chi)(x,\xi)^w}\|_{\mathcal{L}(L^2)}=\|e^{t q(x,\xi)^w}\|_{\mathcal{L}(L^2)}, \ t \geq 0, \num
\end{equation}
because $U$ is a unitary operator on $L^2(\rr^n)$.
Since both operators $iq_2(x'',\xi'')^w$ and $-iq_2(x'',\xi'')^w$ generate contraction semigroups verifying
$$\big(e^{t iq_2(x'',\xi'')^w}\big)^{-1}=e^{t (-iq_2(x'',\xi'')^w)}, \ t \geq 0,$$
the semigroup $e^{itq_2(x'',\xi'')^w}$ is unitary for all $t \geq 0$. It follows from the
tensorization of the variables (\ref{decay1}),
$$e^{t(q \circ \chi)(x,\xi)^w}=e^{tq_1(x',\xi')^w}e^{itq_2(x'',\xi'')^w},$$
and (\ref{decay2}) that
$$\|e^{t q(x,\xi)^w}\|_{\mathcal{L}(L^2)}=\|e^{t q_1(x',\xi')^w}\|_{\mathcal{L}(L^2)}, \ t \geq 0.$$
For proving $(i)$, it is therefore sufficient to prove the exponential decay in time for the norm of the contraction semigroup generated by the operator $q_1(x',\xi')^w$.
We have proved in Proposition~\ref{sperp} (see also (\ref{inf12})) that the spectrum of the operator $q_1(x',\xi')^w$ is only composed of the following eigenvalues
\begin{equation}\label{decay2.5}\inc
\sigma\big{(}q_1(x',\xi')^w\big{)}=\Big\{ \sum_{\substack{\lambda \in \sigma(F_1), \\  \textrm{Re}(-i \lambda) <0} }
{\big{(}r_{\lambda}+2 k_{\lambda}
\big{)}(-i\lambda) : k_{\lambda} \in \nn}
\Big\}, \num
\end{equation}
where $F_1=F|_{S^{\sigma \perp}}$ is the Hamilton map associated to the quadratic form $q_1$ and $r_{\lambda}$
is the dimension of the space of generalized eigenvectors of $F_1$ in $\cc^{2n'}$
belonging to the eigenvalue $\lambda \in \cc$. We have also seen in the proof of Proposition~\ref{sperp} that the contraction semigroup generated by the operator $q_1(x',\xi')^w$
is compact for all $t>0$. We proved this fact in (\ref{kps41}) on the FBI transform side. This allows us to apply
Theorem~2.20 in \cite{davies} to obtain the following description of the spectrum,
$$\sigma(e^{t q_1(x',\xi')^w})=\{0\} \cup \big\{e^{t \mu} : \mu \in \sigma\big(q_1(x',\xi')^w \big)\big\}.$$
Its spectral radius is therefore given by
\begin{equation}\label{decay3}\inc
\textrm{rad}\big(e^{t q_1(x',\xi')^w}\big)=e^{-at }, \num
\end{equation}
where
\begin{equation}\label{decay4}\inc
a=\textrm{inf}\Big\{ \sum_{\substack{\lambda \in \sigma(F_1), \\  \textrm{Re}(-i \lambda) <0} }
{\big{(}r_{\lambda}+2 k_{\lambda}
\big{)}\big(-\textrm{Re}(-i\lambda)\big) : k_{\lambda} \in \nn}
\Big\}. \num
\end{equation}
It follows that the constant $a$ is positive. Since from
Theorem~1.22 in~\cite{davies}, we have
$$-a=\lim_{t \rightarrow +\infty}{\frac{1}{t} \log\|e^{tq_1(x',\xi')^w}\|_{\mathcal{L}(L^2)}},$$
we obtain that there exists $M>0$ such that
\begin{equation}\label{ocello0.1}\inc
\|e^{tq_1(x',\xi')^w}\|_{\mathcal{L}(L^2)} \leq M e^{-\frac{a}{2}t}, \num
\end{equation}
for all $t \geq 0$. This proves that $(iii)$ implies $(i)$. Finally, the fact that $(i)$ implies $(ii)$ is a consequence of a property that we have already mentioned, namely,
when the real part of the symbol $q$ is identically equal to zero then the contraction semigroup
$$e^{tq(x,\xi)^w},$$
is unitary for all $t \geq 0$. This ends the proof of Theorem \ref{theorem3}. $\Box$

\subsection{Spectra of non-elliptic quadratic operators}
\init

In this section, we prove Theorem \ref{theorem2}. Let us consider a complex-valued quadratic form
$$q : \rr_x^n \times \rr_{\xi}^n \rightarrow \cc, \ n \in \nn^*,$$
with a non-positive real part
$$\textrm{Re }q \leq 0,$$
which is elliptic on its singular space $S$,
\begin{equation}\label{til1}\inc
(x,\xi) \in S, \ q(x,\xi)=0 \Rightarrow (x,\xi)=0. \num
\end{equation}
Let us recall that this assumption of partial ellipticity on the singular space ensures that the singular space has a symplectic structure. We can therefore resume the beginning of
our reasoning explained in section~\ref{compact}:

By using the symplectic decomposition of the symbol obtained in section~\ref{sympl}, we deduce from (\ref{sm12}) and (\ref{sm13}) that there exists
$\chi$, a real linear symplectic transformation of $\rr^{2n}$, such that
\begin{equation}\label{til2}\inc
(q \circ \chi)(x,\xi)= q_1(x',\xi')+ iq_2(x'',\xi''), \ (x,\xi)=(x',x'';\xi',\xi'') \in \rr^{2n}, \num
\end{equation}
where $q_1$ is a complex-valued quadratic form on $\rr^{2n'}$ with a non-positive real part
\begin{equation}\label{til3}\inc
\textrm{Re }q_1 \leq 0, \num
\end{equation}
and $q_2$ is a real-valued quadratic form verifying the properties stated in Propositions~\ref{proposition1} and \ref{sperp}.

To obtain the result of Theorem \ref{theorem2}, let us notice from Proposition \ref{proposition1}
that when the symbol $q$ is elliptic on $S$, we can assume that
\begin{equation}\label{til3.1}\inc
q_2(x'',\xi'')=\eps \sum_{j=1}^{n''}{\lambda_j(\xi_j''^2+x_j''^2)},\num
\end{equation}
where $\eps \in \{\pm 1\}$ and $\lambda_j>0$ for all $j=1,...,n''$.

By using again the symplectic invariance of the Weyl quantization given by the theorem~18.5.9 in
\cite{hormander},  we can find a metaplectic operator $U$, which is a unitary transformation on $\lde$ and an automorphism of $\mathcal{S}(\rr^n)$ such that
\begin{equation}\label{til3.15}\inc
(q \circ \chi)(x,\xi)^w=U^{-1} q(x,\xi)^w U. \num
\end{equation}
Since the quadratic form $q_2$ is elliptic on $\rr^{2n''}$, we deduce from the theorem~3.5 in \cite{sjostrand} that the spectrum of the operator $iq_2(x'',\xi'')^w$ is only composed
of eigenvalues with finite multiplicity
\begin{equation}\label{til3.2}\inc
\sigma\big{(}iq_2(x'',\xi'')^w\big{)}=\Big\{ \sum_{\substack{\lambda \in \sigma(iF_2), \\  -i \lambda \in \Sigma(iq_2)\setminus \{0\}} }
{\big{(}r_{\lambda}''+2 k_{\lambda}
\big{)}(-i\lambda) : k_{\lambda} \in \nn}
\Big\},\num
\end{equation}
where $F_2$ is the Hamilton map associated to the quadratic form $q_2$ and $r_{\lambda}''$ is
the dimension of the space of generalized eigenvectors of $i F_2$ in $\cc^{2n''}$
belonging to the eigenvalue $\lambda \in \cc$. We notice from
(\ref{inf11.5}), (\ref{inf12}), (\ref{inf13}) and (\ref{sm14}) that
\begin{equation}\label{til4}\inc
F_1=F|_{S^{\sigma \perp}} \textrm{ and } F_2=\frac{1}{i}F|_S. \num
\end{equation}
Let us notice that if $\lambda$ is an eigenvalue of $F_1$, such that $\textrm{Re}(-i\lambda) \leq 0$, then we necessarily have
$$\textrm{Re}(-i\lambda)<0,$$
because if we had $\textrm{Re}(-i\lambda)=0$, it would imply that the Hamilton map $F_1$ has a real eigenvalue and induce, as we saw in (\ref{venice4}), that the singular space of the
symbol $q_1$ is not reduced to $\{0\}$. However, this singular space is trivial by construction (see $(ii)$ in Proposition~\ref{proposition1}). This proves that
$$\textrm{Re}(-i\lambda)<0.$$
By using now that when the numerical range of a quadratic form $\tilde{q}$ is contained in a closed angular sector
$\Gamma$ with a vertex in $0$ and an aperture strictly less than
$\pi$ then $\lambda$ is an eigenvalue of its Hamilton map $\tilde{F}$ if and only if $-\lambda$ is an eigenvalue of $\tilde{F}$, and
$$-i\lambda \in \Gamma \textrm{ or } i \lambda \in \Gamma,$$
(see section~3 in~\cite{hypoelliptic}),  we obtain from (\ref{inf11.5}), (\ref{inf12}), (\ref{inf13}) and (\ref{sm14}) that
\inc
\begin{align*}\label{til4000}
& \qquad \quad \big\{\lambda \in \cc : \lambda \in \sigma(F), -i\lambda \in \cc_- \cup (\Sigma(q|_S) \setminus \{0\})\big\} \num\\
=& \ \big\{\lambda \in \cc : \lambda \in \sigma(F_1), \textrm{Re}(-i\lambda)<0\big\} \sqcup
\big\{\lambda \in \cc : \lambda \in \sigma(i F_2), -i\lambda \in \Sigma(iq_2) \setminus \{0\}\big\},
\end{align*}
where $F$ is the Hamilton map associated to the quadratic form $q$,
$$\Sigma(q|_S)=\overline{q(S)} \textrm{ and } \cc_-=\{z \in \cc : \textrm{Re }z<0\}.$$

We shall now deduce from the tensorization of the variables (\ref{til2}), (\ref{til3.15}) and (\ref{til4}) that the spectrum of the quadratic differential operator $q(x,\xi)^w$ is only composed of eigenvalues with finite multiplicity
\begin{equation}\label{ocello4}\inc
\sigma\big{(}q(x,\xi)^w\big{)}=\Big\{ \sum_{\substack{\lambda \in \sigma(F), \\  -i \lambda \in \cc_- \cup (\Sigma(q|_S) \setminus \{0\})} }
{\big{(}r_{\lambda}+2 k_{\lambda}
\big{)}(-i\lambda) : k_{\lambda} \in \nn}
\Big\},\num
\end{equation}
where $r_{\lambda}$ is the dimension of the space of generalized eigenvectors of $F$ in $\cc^{2n}$
belonging to the eigenvalue $\lambda \in \cc$. This result is trivial when the singular space is equal to the whole phase space because in that case the quadratic form $q_1$ is identically
equal to 0 and $n''=n$.
We therefore assume in the following that
\begin{equation}\label{ocello4.1}\inc
S \neq \rr^{2n}.\num
\end{equation}
Let us begin by recalling that we know from (\ref{decay2.5}) that the spectrum of the operator
$q_1(x',\xi')^w$ is only composed of eigenvalues with finite multiplicity
\begin{equation}\label{ocello4.3}\inc
\sigma\big{(}q_1(x',\xi')^w\big{)}=\Big\{ \sum_{\substack{\lambda \in \sigma(F_1), \\  \textrm{Re}(-i \lambda)<0} }
{\big{(}r_{\lambda}'+2 k_{\lambda}
\big{)}(-i\lambda) : k_{\lambda} \in \nn}
\Big\},\num
\end{equation}
where $r_{\lambda}'$ is the dimension of the space of generalized eigenvectors of $F_1$ in $\cc^{2n'}$
belonging to the eigenvalue $\lambda \in \cc$. It follows from (\ref{ocello4.3}) that when a positive constant $a$ verifies
\begin{equation}\label{ocello1}\inc
\sigma\big(q_1(x',\xi')^w\big) \cap \{z \in \cc : \textrm{Re }z=-a\} =\emptyset, \num
\end{equation}
the operator $q_1(x',\xi')^w$ only has a finite number of eigenvalues in the half-plane
\begin{equation}\label{ocello2}\inc
\{z \in \cc : -a \leq \textrm{Re }z \}.\num
\end{equation}
In the following, we besides assume that the singular space is not reduced to zero
\begin{equation}\label{ocello4.2}\inc
S \neq \{0\},\num
\end{equation}
because the description (\ref{ocello4}) is then a direct consequence of (\ref{ocello4.3}).
We now need some estimates for the resolvent of the operator $q_1(x',\xi')^w$ to obtain the description (\ref{ocello4}) for the spectrum of the operator $q(x,\xi)^w$.

\bigskip

\begin{proposition}\label{propolol}
For all $a>0$ such that
$$\sigma\big(q_1(x',\xi')^w\big) \cap \{z \in \cc : \emph{\textrm{Re }}z=-a\} =\emptyset,$$
there exists $C_a>0$ such that
\begin{equation}\label{myrtaeq3.1}\inc
\big\|\big(z-q_1(x',\xi')^w\big)^{-1}\big\| \leq C_a, \num
\end{equation}
for all $z \in \cc$ with $-a < \Re z $ and $\abs{\Im z} \geq C_a$. Here the norm is the operator norm on $L^2$.
\end{proposition}

\bigskip

\begin{proof}
When proving Proposition \ref{propolol}, we first recall from (\ref{ocello0.1}) and (\ref{ocello4.1}) that there exist $\tilde{a}>0$ and $M>0$ such that
\begin{equation}\label{ocello5.1}\inc
\norm{e^{tq_1(x',\xi')^w}}\leq M e^{-\tilde{a} t},\ t \geq 0.\num
\end{equation}
By using Theorem 2.8 in \cite{davies}, we can write that for all $z \in \cc$ such that $\textrm{Re }z >-\tilde{a}$,
$$\big(z-q_1(x',\xi')^w\big)^{-1} = \int_0^{+\infty} e^{-z t}\,e^{tq_1(x',\xi')^w} dt,
$$
and we deduce from (\ref{ocello5.1}) that
\begin{equation}\label{myrtaeq3}\inc
\big\|\big(z-q_1(x',\xi')^w\big)^{-1}\big\| \leq \frac{M}{\tilde{a}+\Re z}, \num
\end{equation}
for all $z \in \cc$, $- \tilde{a} < \Re z$. This proves the estimate (\ref{myrtaeq3.1}) when the positive constant $a$ is small enough.
To prove the result in the general case, we shall follow an argument used by L.S. Boulton in~\cite{boulton}. Let us consider a positive constant $a$ verifying (\ref{ocello1}).
We have already seen that the operator $q_1(x',\xi')^w$ has only a finite number of eigenvalues with finite multiplicity in the half-plane
$$\{z \in \cc : -a \leq \textrm{Re }z \}.$$
We can therefore consider $\Pi_a$, the finite-rank spectral projection associated to the eigenvalues
$$\sigma\big(q_1(x',\xi')^w\big) \cap \{ z \in \cc :  -a \leq \Re z \},$$
and write for all $z \in \cc$ with $z \not\in \sigma\big(q_1(x',\xi')^w\big)$ that
\begin{equation}\label{laet5}\inc
\big(z-q_1(x',\xi')^w\big)^{-1} = \big(z-q_1(x',\xi')^w\big)^{-1} \Pi_a + \big(z-q_1(x',\xi')^w\big)^{-1}(1-\Pi_a). \num
\end{equation}
Here
\begin{equation}\label{laet6}\inc
\big(z-q_1(x',\xi')^w\big)^{-1}(1-\Pi_a) = \big(z-q_1(x',\xi')^w|_{{\rm Ran}(1-\Pi_a)}\big)^{-1} (1-\Pi_a), \num
\end{equation}
and we can write by using Theorems 1.22 and 2.8 in \cite{davies} that
\begin{equation}\label{laet2}\inc
\big(z-q_1(x',\xi')^w|_{{\rm Ran}(1-\Pi_a)}\big)^{-1}= \int_0^{+\infty} e^{-z t} e^{tq_1(x',\xi')^w|_{{\rm Ran}(1-\Pi_a)}} dt, \num
\end{equation}
for all $z \in \cc$ with $-b <\textrm{Re }z $ if we set
\begin{equation}\label{laet1}\inc
-b=\lim_{t\rightarrow +\infty} \frac{1}{t} \log \|e^{tq_1(x',\xi')^w|_{{\rm Ran}(1-\Pi_a)}}\|.\num
\end{equation}
Let us now notice that the contraction semigroup
$$e^{tq_1(x',\xi')^w|_{{\rm Ran}(1-\Pi_a)}}=e^{tq_1(x',\xi')^w}(1-\Pi_a),$$
is compact for any $t>0$ since we have seen in section \ref{beha} that the contraction semigroup generated by the operator $q_1(x',\xi')^w$ is compact for any $t>0$.
Since on the other hand the spectrum of the operator $q_1(x',\xi')^w|_{{\rm Ran}(1-\Pi_a)}$ is equal to
$$\sigma\big(q_1(x',\xi')^w\big) \cap \{ z \in \cc : \Re z \leq -a\},$$
we deduce from Theorems 1.22 and 2.20 in \cite{davies} that $-b< -a,$
which implies in view of (\ref{laet1}) that there exists $\tilde{M}>0$ such that
\begin{equation}\label{laet3}\inc
\|e^{tq_1(x',\xi')^w|_{{\rm Ran}(1-\Pi_a)}}\| \leq \tilde{M} e^{-a t},\ t \geq 0.\num
\end{equation}
We then deduce from (\ref{laet2}) and (\ref{laet3}) that
\begin{equation}\label{laet4}\inc
\big\|\big(z-q_1(x',\xi')^w|_{{\rm Ran}(1-\Pi_a)}\big)^{-1}(1-\Pi_a)\big\| \leq \frac{\tilde{M}}{a+\Re z} \|1-\Pi_a\|,\num
\end{equation}
for all $z \in \cc$ with $-a < \textrm{Re }z$.
Since on the other hand
$$\big(z-q_1(x',\xi')^w\big)^{-1} \Pi_a = \big(z-q_1(x',\xi')^w|_{{\rm Ran}\,\Pi_a}\big)^{-1} \Pi_a,$$
and the vector space ${\rm Ran}\,\Pi_a$ is finite-dimensional, we therefore have
\begin{equation}\label{laet7}\inc
\big\|\big(z-q_1(x',\xi')^w\big)^{-1} \Pi_a \big\|=\mathcal{O}_a(1),\num
\end{equation}
for any $z \in \cc$ when $|\Im z|$ is large enough depending on $a$. We finally deduce the result of Proposition \ref{propolol} from (\ref{laet5}), (\ref{laet6}), (\ref{laet4})
and (\ref{laet7}).
\end{proof}

\bigskip

\noindent
\textit{Remark.} Let us notice that the previous proof actually shows that when $q$ is a complex-valued quadratic form on $\rr^{2n}$, $n \geq 1$, with non-positive real part and a zero
singular space
$$S=\{0\},$$
then
$$e^{t q(x,\xi)^w}=e^{t q(x,\xi)^w} \Pi_a+\mathcal{O}_a(e^{-at}), \ t \geq 0,$$
in $\mathcal{L}(L^2)$ for any $a>0$ such that
$$\sigma\big(q(x,\xi)^w\big) \cap \{z \in \cc : \textrm{Re }z=-a\} =\emptyset.$$

\bigskip

We can now resume our proof of Theorem~\ref{theorem2}. In doing so, we recall that any
quadratic differential operator $\tilde{q}(x,\xi)^w$ whose symbol has a non-positive real part, is defined  by the
maximal closed realization on $L^2(\rr^n)$ with the domain
$$\{u \in L^2(\rr^n) :  \tilde{q}(x,\xi)^w u \in L^2(\rr^n)\},$$
which coincides with the graph closure of its restriction to $\mathcal{S}(\rr^n)$,
$$\tilde{q}(x,\xi)^w : \mathcal{S}(\rr^n) \rightarrow \mathcal{S}(\rr^n).$$
By noticing from (\ref{til3.15}) that
\begin{equation}\label{laet9}\inc
\sigma\big((q \circ \chi)(x,\xi)^w\big)=\sigma\big(q(x,\xi)^w\big),\num
\end{equation}
we deduce from (\ref{til2}), (\ref{til3.2}), (\ref{til4000}) and (\ref{ocello4.3}) that it is sufficient
for obtaining (\ref{ocello4}) and ending the proof of Theorem~\ref{theorem2} to establish that
\begin{equation}\label{laet8}\inc
\sigma\big(q(x,\xi)^w\big)=\sigma\big(q_1(x',\xi')^w\big) +\sigma\big(iq_2(x'',\xi'')^w\big).  \num
\end{equation}

We then notice that since the spectra of the operators $q_1(x',\xi')^w$ and $q_2(x'',\xi'')^w$ are only composed of eigenvalues,
we directly get the first inclusion
$$\sigma\big(q_1(x',\xi')^w\big) + \sigma\big(iq_2(x'',\xi'')^w\big) \subset \sigma\big(q(x,\xi)^w\big),$$
by considering the functions
$$u(x',x'')=u_1(x')u_2(x'') \in L^2(\rr^n),$$
where $u_1$ and $u_2$ are respectively some eigenvectors of the operators $q_1(x',\xi')^w$ and $q_2(x'',\xi'')^w$,
since these functions are eigenvectors for the
operator $q(x,\xi)^w$.
We shall now prove the opposite inclusion. Let us consider $z\in \comp$ such that
\begin{equation}\label{ocello113}\inc
z \not\in \sigma\big(q_1(x',\xi')^w\big) +\sigma\big(iq_2(x'',\xi'')^w\big). \num
\end{equation}
In view of (\ref{laet9}), it is sufficient to prove that the map
$$(q \circ \chi)(x,\xi)^w-z : \{u \in L^2(\rr^n) : (q \circ \chi)(x,\xi)^w u \in L^2(\rr^n) \} \rightarrow L^2(\real^n),$$
is bijective to obtain the second inclusion.
We denote by
$$\varphi_{\alpha}(x)= H_{\alpha}(x) e^{-x^2/2}, \ \alpha\in \nat^n,$$
the orthonormal basis of $L^2(\real^n)$ composed by Hermite functions. Here the Hermite polynomials $H_{\alpha}(x)$
satisfy
$$H_{\alpha}(x)=\prod_{j=1}^n H_{\alpha_j}(x_j),$$
and therefore we write
\begin{equation}\label{laeteq8}\inc
\varphi_{\alpha}(x)=\varphi_{\alpha'}(x')\varphi_{\alpha''}(x''),\ \alpha'\in \nat^{n'},\ \alpha''\in \nat^{n''}. \num
\end{equation}
Let us consider the following equation with $u$ and $v$ in $L^2(\real^n)$,
\begin{equation}\label{laeteq9}\inc
(q \circ \chi)(x,\xi)^wu-z u=v. \num
\end{equation}
We can write
\begin{equation}\label{laeteq10}\inc
u(x)=\sum_{\alpha',\alpha''} a_{\alpha'\alpha''}
\varphi_{\alpha'}(x')\varphi_{\alpha''}(x''),\ v(x)=\sum_{\alpha',\alpha''} b_{\alpha'\alpha''}
\varphi_{\alpha'}(x')\varphi_{\alpha''}(x''), \num
\end{equation}
where the two sums are taken for $(\alpha',\alpha'') \in \nn^{n'} \times \nn^{n''}$.
By using from (\ref{til2}) that
$$(q \circ \chi)(x,\xi)^w=q_1(x',\xi')^w+i q_2(x'',\xi'')^w,$$
we obtain from (\ref{til3.1}) that
\inc
\begin{align*}\label{laeteq11}
& \ (q \circ \chi)(x,\xi)^wu-z u \\
=& \ \sum_{\alpha',\alpha''} a_{\alpha'\alpha''}\big[ q_1(x',\xi')^w \varphi_{\alpha'}(x')+ i \mu_{\alpha''}\varphi_{\alpha'}(x')-z \varphi_{\alpha'}(x') \big]\varphi_{\alpha''}(x''),\num
\end{align*}
with
$$\mu_{\alpha''}=\eps \sum_{j=1}^{n''} \lambda_j \left(2\alpha''_j+1\right),$$
since
$$q_2(x'',\xi'')^w \varphi_{\alpha''}(x'')=\mu_{\alpha''}\varphi_{\alpha''}(x'').$$
By setting for any $\alpha''\in \nat^{n''}$,
$$v_{\alpha''}(x')=\sum_{\alpha' \in \nn^{n'}}b_{\alpha'\alpha''}\varphi_{\alpha'}(x') \in L^2(\rr^{n'}),$$
so that according to (\ref{laeteq10}),
\begin{equation}\label{laeteq14}\inc
v(x)=\sum_{\alpha'' \in \nn^{n''}}v_{\alpha''}(x')\varphi_{\alpha''}(x''),\num
\end{equation}
we deduce from (\ref{laeteq11}) that for solving the equation (\ref{laeteq9}), we have to solve all the equations
\begin{equation}\label{laeteq12}\inc
\big(q_1(x',\xi')^w+i\mu_{\alpha''}-z\big)u_{\alpha''}(x')=v_{\alpha''}(x'), \ \alpha'' \in \nat^{n''}, \num
\end{equation}
where
$$u_{\alpha''}(x')=\sum_{\alpha' \in \nn^{n'}} a_{\alpha'\alpha''} \varphi_{\alpha'}(x') \in L^2(\rr^{n'}).$$
We deduce from (\ref{ocello113}) that there is a unique solution $u_{\alpha''}(x')$ in $L^2(\rr^{n'})$ for each of the
equations (\ref{laeteq12}). This proves that for every $v \in L^2(\rr^n)$,
there is at most one solution to the equation (\ref{laeteq9}). Let us denote by $u_{\alpha''}$ the solutions to the equations (\ref{laeteq12}) and
\begin{equation}\label{laeteq13}\inc
u=\sum_{\alpha'' \in \nn^{n''}} u_{\alpha''}(x')\varphi_{\alpha''}(x''). \num
\end{equation}
The equation (\ref{laeteq9}) will have a unique solution in $L^2(\rr^n)$ for every $v \in L^2(\rr^n)$ if we prove that the function $u$ defined in (\ref{laeteq13}) actually
belongs to $L^2(\rr^n)$. This is the case. Indeed, we obtain from (\ref{laeteq14}) and (\ref{laeteq12}) that
\begin{multline*}
\|u\|_{L^2(\rr^n)}^2 = \sum_{\alpha'' \in \nn^{n''}} \|u_{\alpha''}\|_{L^2(\rr^{n'})}^2 =  \sum_{\alpha'' \in \nn^{n''}} \big\|\big(q_1(x',\xi')^w+i\mu_{\alpha''}-z\big)^{-1}v_{\alpha''} \big\|_{L^2(\rr^{n'})}^2 \\
\leq C \sum_{\alpha'' \in \nn^{n''}} \|v_{\alpha''}\|_{L^2(\rr^{n'})}^2=C\|v\|_{L^2(\rr^n)}^2<+\infty,
\end{multline*}
because we deduce from Proposition \ref{propolol} and (\ref{ocello113}) that the quantities
$$\big\|\big(q_1(x',\xi')^w+i\mu_{\alpha''}-z\big)^{-1}\big\|,$$
are bounded with respect to the parameter $\alpha''$ in $\nn^{n''}$. This ends our proof of Theorem \ref{theorem2}.

\bigskip
\bigskip

\noindent
\textsc{UCLA Department of Mathematics, Los Angeles, CA 90095-1555, USA}\\
\textit{E-mail address:} \textbf{hitrik@math.ucla.edu}, \textbf{karel@math.ucla.edu}


\begin{thebibliography}{aa}
\bibitem{boulton}
L.S. Boulton, \textit{Non-self-adjoint harmonic oscillator semigroups and pseudospectra},
J. Operator Theory, \textbf{47}, 413-429 (2002).

\bibitem{davies}
E.B.Davies, \textit{One-Parameter Semigroups}, Academic Press, London (1980).

%\bibitem{DeSjZw} N.Dencker, J.Sj\"ostrand, M.Zworski, \textit {Pseudo-\-spectra of se\-mi\-clas\-si\-cal
%(pseudo)\-diffe\-ren\-tial ope\-rators}, Comm. Pure Appl. Math. {\bf 57}, 384-415 (2004)

%\bibitem{ET}
%M.Embree, N.Trefethen, \textit{Spectra and pseudospectra. The behavior of nonnormal operators and matrices},
%Princeton University Press, 2005.

\bibitem{HN}
B.Helffer, F.Nier, \textit{Hypoelliptic estimates and spectral theory for Fokker-Planck operators and Witten laplacians},
SLN 1862, Springer Verlag, 2005.

\bibitem{HeSjSt} F.H\'erau, J.Sj\"ostrand, C.Stolk, {\it
Semiclassical analysis for the Kramers-Fokker-Planck equation},
Comm. PDE, \textbf{30}, no.4-6, 689-760 (2005).

\bibitem{HeHiSj} F.H\'erau, M.Hitrik, J.Sj\"ostrand, {\it Tunnel effect for Kramers-Fokker-Planck type operators},
Ann. Henri Poincar\'e, to appear.

\bibitem{hypoelliptic}
L.H\"{o}rmander, \textit{A class of hypoelliptic pseudodifferential operators with double characteristics},
Math. Ann., \textbf{217},
165-188 (1975).

\bibitem{hormander}
L.H\"{o}rmander, \textit{The analysis of linear partial differential operators} (vol. I,II,III,IV),
Springer Verlag (1985).

\bibitem{mehler}
L.H\"{o}rmander, \textit{Symplectic classification of quadratic forms, and general Mehler formulas},
Math. Z., \textbf{219}, 413-449 (1995).

\bibitem{MelinSj} A.Melin, J.Sj\"ostrand, {\it Determinants of
pseudodifferential operators and complex deformations of phase space},
Methods Appl. Anal., \textbf{9}, no. 2, 177-237 (2002).

\bibitem{mz}
K.Pravda-Starov, \textit{Contraction semigroups of elliptic quadratic differential operators},
to appear in Mathematische Zeitschrift (2007).

%\bibitem{karel}
%K.Pravda-Starov, \textit{On the pseudospectrum of elliptic quadratic differential operators}, preprint, 2007.

\bibitem{sjostrand}
J.Sj\"{o}strand, \textit{Parametrices for pseudodifferential operators with multiple characteristics},
Ark. f\"{o}r Mat., \textbf{12}, 85-130 (1974).
\bibitem{Sj82} J.Sj\"ostrand, {\it Singularit\'es analytiques microlocales}, Astérisque, \textbf{95}, 1-166,
Soc. Math. France, Paris (1982).

\bibitem{SjHokk} J.Sj\"ostrand, {\it Analytic wavefront sets and operators with multiple characteristics}, Hokkaido
Math. Journal, {\bf 12}, no. 3, part 2, 392-433 (1983).

\bibitem{Sj95} J.Sj\"ostrand, {\it Function spaces associated to
global I-Lagrangian manifolds}, Structure of solutions of differential
equations, Katata/Kyoto, 1995, World Sci. Publ., River Edge, NJ (1996).

\bibitem{SjZw} J.Sj\"ostrand, M.Zworski, {\it The complex scaling method for scattering by strictly convex obstacles},
Ark. f\"or Mat., \textbf{33}, 135-172 (1995).

\end{thebibliography}
\end{document}